# Multi-agent interaction in the problem of territorial distribution of production points and temporary storage facilities: the model's construction and analysis.

Oleg Malafeyev[1], Schenikova Snezhana[2]

[a]Saint-Petersburg State University, Russia

**Abstract**

This article discusses the algorithms for finding the optimal solution of problems related to the location of temporary storage of goods, warehouses, factories for processing raw materials and shops selling the final product in the transport network. An algorithm is also proposed for finding a compromise solution to the problem of maximizing profits for each agent.

## 1. Literature Review

To date, there are many studies on the construction of models for multi-agent interaction. For example Umar Manzoor, Maria Zubair and Kanwal Batool[61] studied an efficient distribution of resources such as medicine, food, water etc. They have proposed a solution in which a big disaster region is divided into smaller areas (regions) and multi-agent system are used to deliver food in these regions. The idea is that agents use the shortest path algorithms and coordination for efficient food delivery within their region. Another work of Maria Caridi and Sergio Cavalieri[62] is dedicated to the impact of multi-agent systems on manufacturing practices in the enterprise and at a broader level of the supply chain. The Michal Pechoucek, Ales Riha, Jiri Vokrinek, Vladimir Marik and Vojtech Prazma[63] consider such a multi-agent production planning technology which makes it easier to optimize the use of resources and the supply chain while meeting customer needs. The J. Li, J.Y.H. Fuh, Y.F. Zhang and Andrew Y C Nee[64] investigated the design and production planning system based on several agents. As a result, a prototype of a distributed collaborative design environment was proposed. The J. Rouchier, F. Bousquet, O.

[1] malafeyevoa@mail.ru
[2] snezha_00_97@mail.ru

Barreteau, C. Le Pagez and J.-L. Bonnefoy in the article[65] mentioned dependence of the distribution of goods between agents on the distribution of agents in space. As a result, four societies were built using multi-agent simulation models that address issues related to the use of conventional renewable goods

# 2. Introduction

We would like to consider the case of a single warehouse for several small or medium-sized firms, enterprises. Usually this situation arises when sales are limited to one or several close regions. Usually there are no problems here. Questions arise when working with large firms that have a large international market. What warehouses in this case are selected in the optimal way? To solve this problem, the compromise help method is used. The study of warehouse space needs in different sales regions is also significant. There are two types of placement of the warehouse network: centralized and decentralized. In the first case, the presence of one large warehouse, in the second - several warehouses, scattered in different sales regions.

The capacity of material flows, their rational organization, demand in the sales market, the size of the sales region, the concentration of consumers in it, the location of suppliers and customers relative to each other are the main criteria by which the location of warehouses is determined, their number. From this arises another mathematical problem in the formation and territorial location of the warehouse network in an optimal way. You can build a new warehouse or buy an existing warehouse and operate it. It requires significant capital investment. On the other hand, due to the maximum approach of warehouses to customers, there is a reduction in distribution costs. This implies the problem of optimal placement of temporary storage of products and raw materials. It is also necessary to solve the problem of the territorial location of factories for the processing of raw materials and shops where products are sold. Thus, there are several optimization problems for which you can find a solution using the methods proposed by us. The ideas and approaches from [5][40][56] are used in the paper.

# 3. Informal statement of the problem

Suppose there is a certain plane $\boldsymbol{\pi}$ on which some finite transport network

**C** = $(N, \mathrm{p}, \mathrm{d})$, where $\boldsymbol{N} = \{x_1, x_2, \dots, x_k\}$ - a finite set of vertices, $\boldsymbol{p}$ - a function of transport costs, $\boldsymbol{d}$ - limited capacities. Denote by $\boldsymbol{d} = \boldsymbol{d}\,(d_1,..,d_{N_A})$ and $\boldsymbol{d'} = \boldsymbol{d'}(d'_1\,,..,d'_{N_B})$ capacity for each type of raw materials and manufactured products. $N_A$, $N_B$ is the number of types for each type of raw material respectively. Similarly, we denote the functions of transport costs: $\boldsymbol{p}$ = $\boldsymbol{p}(p_1,..,p_{N_A})$ and $\boldsymbol{p}'$=$\boldsymbol{p}'(p_1'$ ,.., $p_{N_B}')$. Consider the following scheme: "the extractionof raw materials → its storage point → production → point of temporary storage of products → store". In this regard, the edge can not be present at once two parameters $\boldsymbol{d}$ and $\boldsymbol{d}'$ or $\boldsymbol{p}$ and $\boldsymbol{p}'$respectively. Further we will designate it only for $\boldsymbol{p}$ and $\boldsymbol{d}$, which have a different meaning depending on the context. We will also introduce such concepts as points of production (**A** = $\{A_1,\ldots,A_{N_A}\}$) and raw materials storage points (**S** = =$\{S_1,\ldots,S_{k_A}\}$), raw material processing plants (**B**=$\{B_1,\ldots,B_{N_B}\}$) and stores that sell manufactured goods (**M**=$\{M_1,\ldots,M_m\}$).

Further, we will assume that our points of extraction of raw materials are fixed, that is, they are located at strictly defined nodes. For each plant, a set of different types of raw materials is set, from which one type of manufactured products is obtained. Every buyer wants to buy some set of goods. Goods are delivered from the warehouse (to factories and shops), which corresponds to a given type of raw materials and goods. Also introduced the concept of value of each manufactured product. It represents the sum of the value of the goods at the point of extraction of raw materials, transportation costs for its delivery first to the storage warehouse, and then to the factory, production costs, costs for delivering goods to the store and storage costs in warehouses.

Several agents are being considered. They own production facilities, transportation, road junctions and raw materials extraction points. Also important is the concept of action. Each agent receives a certain percentage of the profits in each area. If there is no stock, we deal with zero income. Determine what the costs are for each type of agent. Those involved in road interchanges and transportation spends money on road maintenance, as well as servicing or replacing vehicles and wages for workers. Those who own production facilities and points of extraction of raw materials, deal with some fixed costs: payment for gas, electricity, water, maintenance costs, rent, and salary. Those who deal with warehouses spend money on maintenance (insurance, fire safety, wages to workers, and maintenance of equipment in a warehouse). The size of the variable cost is small, so it can be omitted.

After the introduction of basic concepts and notation, we define a mathematical problem. It represents a mathematical problem of minimizing costs for each of the agents. It can be seen that this problem is a multi-agent interaction problem. We remember that the territorial location of production points and food storage points is also part of our mathematical problem. We introduce some parameter of the total profit and subtract from it all the costs we considered. Thus, our minimization problem is reduced to its dual problem - the problem of maximizing profits for each agent. The problem is divided into the following stages:

1. Compiling a mathematical description for the given conditions;

2. Introduction of cost functions (for production, transport, production);

3. Description of algorithms and their application in determining the optimal solutions and equilibrium positions for a given problem;

## 4. Formalization, construction and analysis of the multi-agent interaction model

Let a finite transport network (*N*, *p*, *d*) be given on a plane, where ***N*** is a finite set of nodes, ***p*** are the transport cost functions (***p***: (***N***, ***N***) → $R^1$), set on network edges, ***d*** − limited throughput (***d***: (***N***, ***N***) → $R^1$). The edge of the network is an ordered pair of nodes ($x_i$, $x_j$). The network is represented as follows:

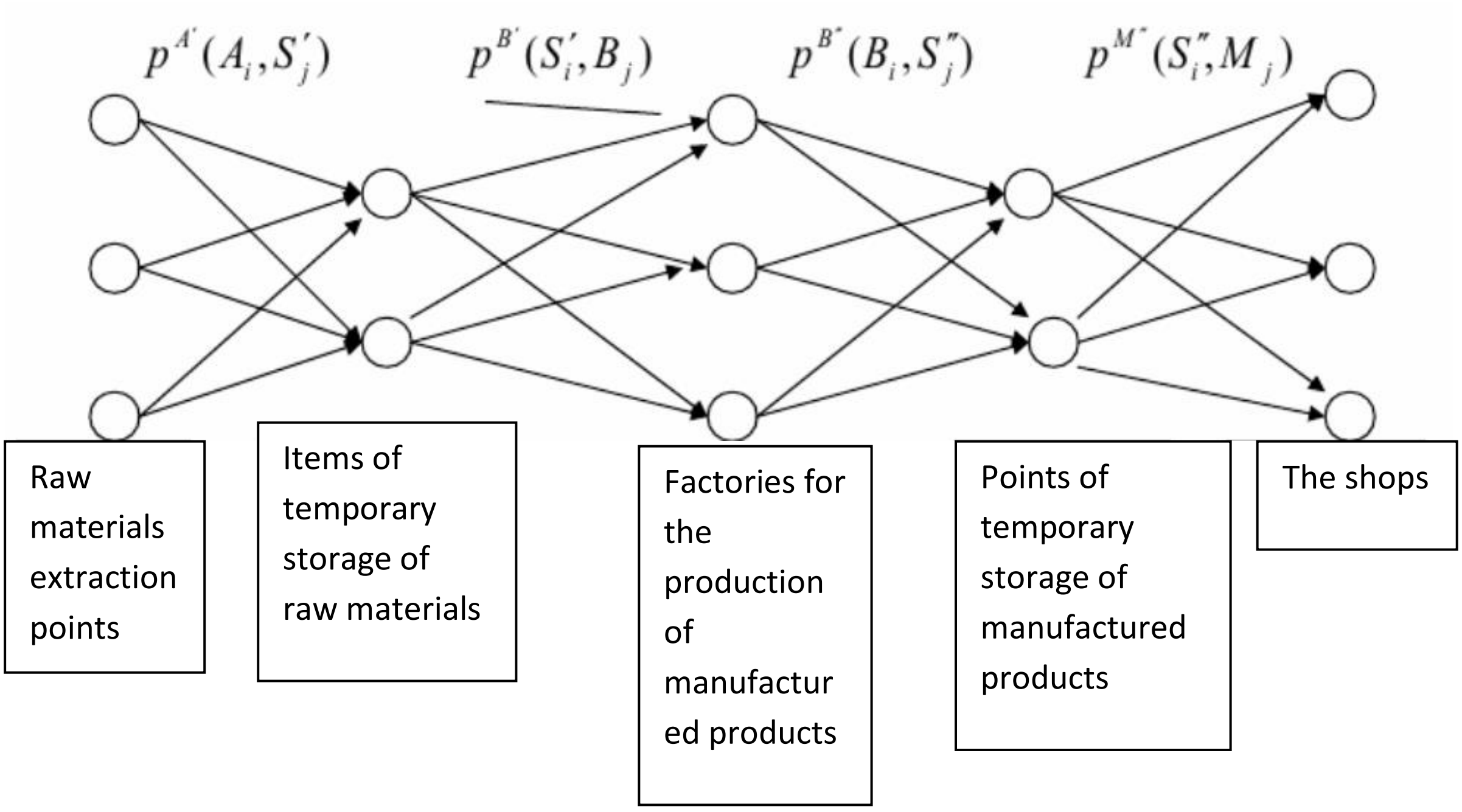

$\boldsymbol{N_A}$ is a number of points of extraction of raw materials, $A_1,\ldots,A_{N_A}$ are the points of extraction of raw materials, $C_i^A = C_i^A(c_1^A,\ldots,c_{N_A}^A)$ is the cost of the i-th set of raw materials, $S_i = S_i(s_1,\dots,s_{k_A})$ is cost function of the i-th set of raw materials. Each edge ($x_i$, $x_j$) of the network between the points of extraction of raw materials ($A_1,\dots,A_{N_A}$) and points of temporary storage of products ($S_1,\dots,S_{k_A}$)) the transportation cost function of moving the i-th unit of raw materials ($p_i^A(A_j, S_k)$, where j=1,…,$N_A$ – number of types of raw materials, k=1,,…,$K_A$ – number of temporary storage of raw materials) matches a non-negative number $p_{ij}^A(A_j, S_k)$ for each type of raw materials. Thus, it is possible to make a matrix $P=(p_{ij})$, i=1,…,$N_B$, j=1,…,$N_A$.

The network has free nodes. They can be located $N_B$ factories for the processing of raw materials into manufactured products $B_1 \dots B_{N_B}$. For the production of a unit of manufactured products, a certain set of raw materials is required $f_i = (V_1,\dots,V_{N_A})$, i=1,…$N_B$ is the number of types of manufactured products. Cost per unit of $C_{i_B}$ production, cost of production of the i-th set of manufactured products: $C_{i_B} = C_{i_B}(C_{1_B},\dots,C_{B_{N_B}})$. Similarly, we denote the functions of transportation costs: $p_i^A(S_k, B_j)$, where j = 1,…,$N_B$ is the number of production points, k = 1, ..., $K_A$ is the number of temporary storage of raw materials. This function assigns to each edge ($x_i$, $x_j$) a nonnegative number $p_{ij}^A(S_k, B_j)$ for each type of raw material.

Further, by analogy, the following notation is introduced:

$p_i^A(S_k, B_j)$ – the function of transport costs between the temporary storage of raw materials and plants for each edge ($x_i$, $x_j$);

$p_{ij}^A(S_k, B_j)$ – function values for each type of raw material, where j = 1, ..., $N_B$ is the number of production points, k = 1 ... $K_A$ is the number of temporary storage of raw materials;

$B_1, \dots, B_{N_B}$ – plants for the production of manufactured products;

$S_i' = S_i'(s_1, \dots, s_{k_A})$ – the cost function of storing the i-th set (of a certain type) of the final raw product;

$p_i^B(B_J, S_k')$ – functions of transportation costs for moving a single i-th product, where j=1…$N_B$ is the number of production points, , k =1…$K_B$ is the number of temporary storage of manufactured products;

$p_{ij}^B(B_J, S_k')$ – function that associates with each edge the transport costs of a certain type of manufactured product;

$M_1, \dots, M_m$ – network nodes in which there are M stores. In each of them defined demand;

$m_i = (w_1, \dots, w_{N_B})$ – a set of manufactured products that are purchased for the i-th store;

$p_j^B(S_k', M)$ – functions of transportation costs for moving products between temporary storage of final products and stores;

$p_1^{M'}(x_i, x_j), \dots, p_M^{M'}(x_i, x_j)$ – the value of the function for each type of manufactured product;

We deal with four items to obtain the costs of storage, transportation and extraction of raw materials:

1. $C^A f_i$ or $\sum_{i=1}^{N_A} C_i^A v_i$ — cost of production of a set of raw materials;

2. $S' f_i$ or $\sum_{i=1}^{N_A} S_i' v_i$— storage costs;

3. $p^A(A_i, S_i')$ or $\sum_{i=1}^{N_A} p_i^A (A_i, S_j') V_i$ – transportation costs for the transportation of raw materials from the extraction point to the storage of products;

4. $p^A(S_i', B_j) f_i$ or $\sum_{i=1}^{N_A} p_i^A (S_i', B_j) V_i$ – transportation costs for the transportation of a set of raw materials from the storage warehouse to the factory;

So, the cost of purchasing a set of raw materials $f_i = (V_1 \dots V_{N_A})$ is:

$$C_k^z(f_i) = C^A f_i + P^A(A_i, B_j) f_i + S f_i ,$$

where k = 1,...,$N_B$ is the number of points of raw materials processing into manufactured products or

$$C_k^z(f_i) = \sum_{j=1}^{N_A} C_j^A v_j + \sum_{j=1}^{N_A} p_j^A (A_l, \ B_m) v_j + \sum_{j=1}^{N_A} S_j' v_j$$

where k = 1,...,$N_B$ is the number of points of raw materials processing into manufactured products.

By analogy, we make the purchase price for the manufactured product:

1. The cost of producing a specific set of manufactured products

$$C_{m_i}^B \; or \; \sum_{i=1}^{N_A} c_i^B \, v_i$$

2. Transportation costs for the transportation of a specific set of manufactured products

$$p^B(B_i, S_j'')m_i \;\; or \; \sum_{i=1}^{N_A} p_i^B \, (B_i, S_j'')w_i$$

3. The cost of storing a specific set of manufactured products

$$S''m_i \; or \sum_{i=1}^{N_A} S_i' \, w_i$$

4. Transportation costs for the transportation of a specific set of manufactured products

$$p_i^B(B_j, S_j'')m_i \; or \sum_{i=1}^{N_A} p_i^B \, (B_i, S_j'')w_i$$

We get the costs to determine the manufactured products for a particular store:

$$C_k^z(f_i) = C^B m_i + p^B(B_{i,} M_j)w_i + \; S''m_i \, ,$$

where k = 1,...,M  is the number of stores where the sale of final products. The total cost of the manufactured product is represented as C = $(C_1^z, \ldots, C_M^z)$. Next, we need to determine the distance between objects.  For this, the concept of the Euclidean norm is introduced:

$$r(x_j, y_j) = \sqrt{\left(x_j - x_k\right)^2 + \left(y_j - y_k\right)^2},$$

where $(x_j, y_j)$ – coordinates of point $x_i$

Let $R_i$ be the admissible distance between objects, i=1,…,$N^2$, $N^2$ – the possible number of edges. Based on all of the above, we can set the minimization problem of the total costs:

$$\eta = CM \rightarrow min$$

It is necessary that the following conditions be met:

1. $r_i \leq R_i$ – distance between objects (warehouse and factory, for example);

2. $f_i \leq di$ – edge capacity;

3. $p_i \geq 0$, i=1,…,$N^2$, where $N^2$ – possible number of network edges – shipping costs;

4. $v_i \geq 0$, i=1,…,$N_A$ – raw material;

5. $\sum_{i=1}^{N_A} v_i \leq V$ – total row material;

6. $w_i \geq 0$, i=1,…,$N_B$ – manufactured product;

7. $\sum_{i=1}^{N_B} w_i \leq W$ – total manufactured product;

We also have three main agents: *Agent A* – an agent owning all production facilities, *Agent B* – an agent owning interchanges and warehouses, *Agent C* – an agent owning stores.

We can make a small conclusion. The problem is that with the concepts, indicators and conditions given above, we need to find a compromise solution that will maximize the profit of each agent.

# 5. Formalization of the algorithm for solving the problem by the example of a particular case

Suppose we have 2 points of raw material extraction – A and B, 2 points of temporary storage of raw materials - C and D, production points – E and F, points of temporary storage of manufactured products – G and I and 4 stores – $M_1,M_2,M_3,M_4$.

## 5.1. Definition of the shortest path matrix between all pairs of points

To begin with, we introduce such a concept as the shortest-path matrix from any node $(x_i)$ of the network graph to any other node:

| | $x_1$ | $x_2$ | ... | $x_n$ |
|---|---|---|---|---|
| $x_1$ | 0 | $a_{x_1x_2}$ | ... | $a_{x_1x_n}$ |
| $x_2$ | $a_{x_2x_1}$ | 0 | ... | $a_{x_nx_2}$ |
| $\vdots$ | $\vdots$ | $\vdots$ | $\ddots$ | $\vdots$ |
| $x_n$ | $a_{x_nx_1}$ | $a_{x_nx_2}$ | ... | 0 |

where $a_{x_ix_j}$ is the path length between nodes $x_i$ and $x_j$. It is clear that the distance between the same nodes is zero.

To obtain this matrix, we will use the Floyd algorithm. This algorithm is useful in our case, since we are dealing with a large number of edge pairs between pairs of vertices. Next, a matrix of total costs is built for each set of products at each point of sale:

| | $x_1$ | $x_2$ | ... | $x_n$ |
|---|---|---|---|---|
| $k_1$ | $k_1(x_1)$ | $k_1(x_2)$ | ... | $k_1(x_n)$ |
| $k_2$ | $k_2(x_1)$ | $k_2(x_2)$ | ... | $k_2(x_n)$ |
| $\vdots$ | $\vdots$ | $\vdots$ | $\ddots$ | $\vdots$ |
| $k_n$ | $k_n(x_1)$ | $k_n(x_2)$ | ... | $k_n(x_n)$ |

where $x_i$ is the i-th shop, $k_i$ is a set of goods (for a specific type), $k_i\ (x_j)$ are total costs, i = 1 ... n, j = 1 ... m. In our case, n = m = 2.

## 5.2. Determination of the functions of winning players

The winning functions of each player are the necessary values for finding a compromise solution in the placement problem. Recall that in our mathematical problem, the players are certain agents. Each of them has production capacity (agent $A_1$), transportation, transport interchanges (agent $A_2$) and raw materials extraction points (agent $A_3$). We formulate an algorithm to search the payoff functions:

**Step 1. Calculate total demand.**

This step allows you to determine the required amount of extracted resources and the amount of production of manufactured products.

**Step 2.Preliminary estimate of revenues and costs at each network node.**

Let the first agent $A_1$ own the points of temporary storage of raw materials, points of storage of final products, transport. The second agent $A_2$ own production points. And the third $A_3$ own shops located in the points of consumption.

Then the profit of the first agent ($P_1$) is calculated by the formula:

$$P_1 = I_1 - C_1,$$

where $I_1$ is the income from providing storage services, $C_1$ is the cost of transporting and storing raw materials and manufactured products

Profit of the second agent ($P_2$):

$$P_2 = I_2 - C_2,$$

where $I_2$ is the income from the cost of manufactured products, $C_2$ is the cost of purchasing processed raw materials.

The output is calculated using the Cobb-Douglas function:

$$Q = JK^{\alpha}L^{\beta},$$

where $J$ is a coefficient of manufactured products, $K^{\alpha}$, $L^{\beta}$ are the costs of raw materials of the form $\alpha$ and $\beta$, respectively.

The profit of the third agent ($P_3$) is calculated as follows:

$$P_3 = I_3 - C_3,$$

where $I_3$ is the income from the sale of final products, $C_3$ is the cost of purchasing goods from the manufacturer.

**Step 3.Consideration of possible options for the location of production points.**

At this step, the net profit of each agent is calculated separately from each other.

## 5.3.Finding a compromise solution. Algorithm Description

We will find a compromise solution. We make the next algorithm:

**Step 1. Construction of the payoff matrix (Γ)**

$$\boldsymbol{\Gamma} = (\alpha_{l,m}),$$

where l is the number of players, m is the number of situations in the game

**Step 2. Compilation of a vector called "ideal" (M)**

$$M = \begin{pmatrix} M_1 \\ \vdots \\ M_l \end{pmatrix}, M_l = \max_m \alpha_{l,m}$$

It consists of the maximum income that producers receive.

**Step 3. Calculation of residuals**

Under the discrepancies we will understand the deviations of the income of each manufacturer from the maximum income:

$$\Gamma_M = \left(M - \alpha_{l,m}\right) = (\beta_{l,m})$$

**Step 4.Arrange income Values**

In each situation, we arrange the income in ascending order so that the first line contains the smallest values, and the last - the largest:

$$\max_m \beta_{m,l} = \max_l (M - \alpha_{l,m})$$

**Step 5. The principle of "minmax"**

Among the found maximum discrepancies choose the minimum value:

$$\min_m \max_l (M - \alpha_{l,m})$$

There are cases when there are several situations in the last line with the same minimum. Then you need to go to the line above and look for the minimum value there. The resulting situations are the desired compromise set (solution).

# 6. Objectives of finding the costs and profits of each individual object in the network. Description.

## 6.1.Production Planning Challenge

We reformulate the conditions of the problem for the points of processing raw materials and manufacturing products. Suppose we have resources of a specific type. We need to produce from them manufactured products, and a certain set of different types. Then they go on sale. We know the demand for goods. According to it, we can make a plan for the extraction and production of goods. Imagine our problem in a mathematical form.

Let ***n*** be the number of product types that the company produces. Then for $U_1,...,U_n$ we denote the types of products themselves. For each type of manufactured product, the company has a plan. According to it, it is obliged to fulfill at least $b_1$ units of final production $U_1$, at least $b_2$ units of final production $U_2$, etc. A case of over-fulfillment of the plan is possible, but even here there are limits. No more than $\beta_1,...,\beta_n$ units of each type must be produced. Let $m$ be the amount of resources extracted. The resources themselves are denoted by $s_1,...,s_m$. Each type of resource extracted is limited by the numbers $\gamma_1,...,\gamma_m$, respectively. For $a_{ij}$ we denote the number of units of the extracted resource of the form $s_i$, i = 1, ..., m, which is necessary for the manufacture of one unit of final production $U_j$, j = 1,...,n. Then you can make the following matrix:

| Extracted resource | Raw materials | | | | | |
|---|---|---|---|---|---|---|
| | $U_1$ | $U_2$ | ... | $U_j$ | ... | $U_n$ |
| $S_1$ | $a_{11}$ | $a_{12}$ | … | $a_{1j}$ | … | $a_{1n}$ |
| $S_2$ | $a_{21}$ | $a_{22}$ | … | $a_{2j}$ | … | $a_{2n}$ |
| ... | … | … | … | … | … | … |
| $S_i$ | $a_{i1}$ | $a_{i2}$ | … | $a_{ij}$ | … | $a_{in}$ |
| … | … | … | … | … | … | … |
| $S_m$ | $a_{m1}$ | $a_{m2}$ | … | $a_{mj}$ | … | $a_{mn}$ |

When implementing the plan, the unit of final production $U_i$ brings the enterprise a profit $c_i$, i = 1, ..., n. It is necessary to plan the production of the final product (its quantity) in an optimal way. Optimal in the sense that the plan must either be completed or exceeded (to the extent permitted). Total profit should go to the maximum. It can be seen that this problem can be represented as a linear programming problem:

- $x_i \geq b_i$ – obligatory execution of the plan, i = 1,...,n;
- $x_i \leq \beta_i$ – absence of excessive production, i = 1, ..., n;
- $xA \leq \gamma$, where $A = \{a_{ij}\}$, i = 1,.., n, j = 1,..,m is the matrix of the number of units finite production, $\gamma = (\gamma_1,.., \gamma_m)$ – limiting reserves of extracted resources, $x = (x_1,..,x_n)$ – the number of units of final production $U_1,..,U_n$, which we will produce -restrictions so that we have enough extracted resources;
- $L = cx \rightarrow max$, where $c = (c_1,..,c_n)$ is the profit from the produced final product, $x = (x_1,..,x_n)$ - the number of units of the final product $U_1,..,U_n$, which we produce - the profit brought by the plan $(x_1,..,x_n)$.

Thus, the production planning problem is formulated. Next, we consider the problem of selective control of manufactured products.

Products must be of high quality. For this, a sampling control system is organized. It is necessary to organize the control in such a way that at the minimum cost of control to ensure the specified level of quality. Let $R$ be the average expected cost of control per unit of time. Then we will consider these costs as a natural indicator of efficiency with the condition that the control system provides the specified level of quality. As a quality level, you can take the average percentage of rejection higher than the specified, i.e. $R \rightarrow min$

## 6.2. The problem of loading temporary storage points

When planning points for temporary storage of products, you should consider the area occupied by the room where the product is stored. Also required transport platform for loading/unloading products. How do you plan to fill the warehouse? Consider further the problem of loading.

This problem is understood as the problem of rational loading of an item where goods are temporarily stored. This item has limitations on volume or carrying capacity. A load placed in a temporary storage facility makes a profit. We need to load the place with such goods, which in the end will bring the total profit. It is clear that this indicator should strive to the maximum.

Let $P$ be a temporary storage facility, $W$ be the volume of stored items, $n$ be the number of product names. Then $m_i$ is the number of products of the i-th name to be loaded. For $r_i$ we denote the profit obtained from one loaded item of the i-th name, $w_i$ is the weight of one item of the i-th name. Again we deal with the problem of linear programming:

Objective function: $z = \sum_{i=1}^{n} r_i \, m_i \rightarrow max$

Conditions: $\sum_{i=1}^{n} w_i \, m_i \leq W, m_1, \ldots, m_n \geq 0;\ m_1, \ldots, m_n \in Z$

To solve this problem, a dynamic programming model is built. Stage $i$ is assigned to the subject of the i-th name, i = 1,...,n

At the i-th stage, the solutions to the problem are described by the number of mi items of the i-th name to be loaded. $r_i m_i$ – the resulting profit. The lower bound of the mi values is 0, the upper bound is $\left[\frac{W}{w_i}\right]$

Let $x_i$ be the total weight of objects at the i-th stage, the loading decisions of which were taken at subsequent stages i, i = 1,...,n. It is the weight restriction that is the only restriction that unites all n stages.

Further, for $f_i(x_i)$, we denote the maximum total profit from stages i, i = 1, ..., n with the given state $x_i$ .

To determine the recurrent equation, the following procedures are introduced:

**Step 1. The expression of the function $f_i(x_i)$ through $f_{i+1}(x_{i+1})$**

$$f_i(x_i) = \max_{\substack{m_i = 0, \ldots, \left[\frac{W}{w_i}\right] \\ x_i = 0, \ldots, W}} (r_i m_i + f_{i+1}(x_{i+1}),\ \ i = 1, \ldots, n,\ where\ \ f_{n+1}(x_{n+1}) \equiv 0$$

**Step 2. The expression $x_{i+1}$ through $x_i$**

By definition, the quantity ($x_i - x_{i+1}$) is the weight loaded at the i-th stage. In other words $x_i - x_{i+1} = w_i m_i$. From here we get $x_{i+1} = x_i - w_i m_i$.

Thus, the recurrence equation has the form:

$$f_i(x_i) = \max_{\substack{m_i = 0, \ldots, \left[\frac{W}{w_i}\right] \\ x_i = 0, \ldots, W}} (r_i m_i + f_{i+1}(x_i - w_i m_i),\ \ i = 1, \ldots n.$$

The problem of resource allocation, in which a particular resource is distributed among a finite number of activities, refers to a stated load problem. The optimality is taken to maximize the profit function. It can be seen that, similar to the problem of loading, the state at the i-th stage is the total amount of the resource, which is distributed at stages i, i + 1,...,n.

## 6.3 Transport mathematical problem

The problem of transportation products – one of the main issues of the modern economy. Various factors affecting production in different places related to location, quality of sources of raw materials cause transportation of products from one place to another. The lack of vehicles, their large load – the reasons for the need to address issues related to transport.

The transport problem is modeled as follows. Let $A_1,\ldots A_m$ – points of departure. The stocks of some homogeneous goods are concentrated in them, $a_1,\ldots a_m$ is the amount of these goods (units). $B_1,\ldots,B_n$ – destinations that require cargo $b_1,\ldots b_n$. An important condition is the balance:

$$\sum_{i=1}^{m} a_i = \sum_{j=1}^{n} b_j$$

Denote by $c_{ij}$ the cost of transportation of cargo from the point of departure of $A_i$ to the destination of $B_j$ (i = 1,...,m; j = 1,...,n). We can create a matrix, all values of which are set in advance:

| | | | |
|---|---|---|---|
| $c_{11}$ | $c_{12}$ | … | $c_{1n}$ |
| $c_{21}$ | $c_{22}$ | … | $c_{2n}$ |
| … | … | … | … |
| $c_{m1}$ | $c_{m2}$ | … | $c_{mn}$ |

We assume that the cost of transporting several units of cargo is proportional to their quantity. The total cost of all shipments should be minimal. Thus, we deal with the minimization problem. It is necessary to make an optimal transportation plan. Namely, from where, where, how many units of cargo to carry to fulfill all applications.

Let $x_{ij}$ be the number of units of cargo that is sent from point $A_i$ to point $B_j$. It is clear that these variables can not take negative values.

Create a matrix of these variables:

| | | | |
|---|---|---|---|
| $x_{11}$ | $x_{12}$ | … | $x_{1n}$ |
| $x_{21}$ | $x_{22}$ | … | $x_{2n}$ |
| … | … | … | … |
| $x_{m1}$ | $x_{m2}$ | … | $x_{mn}$ |

Denote $x_{ij}$ as “carriage” and $c_{ij}$ as “transportation plan” (i = 1,...,m; j = 1,...,n). These variables must satisfy the conditions:

1. The total amount of cargo in the points of departure is equal to the stock of cargo in this paragraph.

Imagine these conditions as a system

$$\begin{cases} x_{11} + x_{12} + \cdots + x_{1n} = a_1 \\ x_{21} + x_{22} + \cdots + x_{2n} = a_2 \\ \ldots \\ x_{m1} + x_{m2} + \cdots + x_{mn} = a_m \end{cases}$$

or

$$\sum_{j=1}^{m} x_{ij} = a_i, \;\; i = 1, \ldots n$$

2. The total amount of cargo delivered to each destination from all points of departure must be equal to the application filed from this point:

$$\begin{cases} x_{11} + x_{12} + \cdots + x_{1n} = b_1 \\ x_{21} + x_{22} + \cdots + x_{2n} = b_2 \\ \ldots \\ x_{m1} + x_{m2} + \cdots + x_{mn} = b_m \end{cases}$$

or

$$\sum_{i=1}^{n} x_{ij} = b_j, \;\; j = 1, \ldots, m$$

And we denote the minimization of the total cost of all shipments for:

$$L = \sum_{i=1}^{m} \sum_{j=1}^{n} c_{ij} x_{ij} \to min$$

Thus, we will collect all the conditions into a single whole and obtain a transport problem — a linear programming problem:

$$\sum_{j=1}^{m} x_{ij} = a_i \,, i = 1, \dots, n, \qquad (1)$$

$$\sum_{i=1}^{n} x_{ij} = b_j \,, j = 1, \dots, m, \qquad (2)$$

$$L = \sum_{i=1}^{m} \sum_{j=1}^{n} c_{ij} x_{ij} \to min \,, \qquad (3)$$

We take into account that the number of basis variables is (m + n-1). This is due to the fact that conditions (1) and (2) are not linearly independent. This follows from the condition:

$$\sum_{i=1}^{m} a_i = \sum_{j=1}^{n} b_j$$

A transportation plan is said to be admissible if it satisfies conditions (1), (2). The supporting admissible plan is such a plan that there are no more than (m + n-1) basic variables in it distinct from zero. The rest of the carriage is zero.

Under the optimal plan $x_{ij}$ (i = 1,...,m; j = 1,...,n), we mean the plan, which among all other plans leads to $L_{min}$ (the minimum total cost of transportation).

There is also a variant of applications – “transportation problem with the wrong balance”. Here the inequality holds:

$$\sum_{i=1}^{m} a_i \neq \sum_{j=1}^{n} b_j$$

Two cases are possible:

1. $\sum_{i=1}^{m} a_i > \sum_{j=1}^{n} b_j$

- the amount of inventory is greater than the amount of needs. That is, we still have some resources. The problem is reduced to the transport problem by entering some "fictitious" destination $B_f$.

It can be seen that this mathematical problem really comes down to a balanced the mathematical problem, since

$$2. \sum_{i=1}^{m} a_i < \sum_{j=1}^{n} b_j$$

that is, there are not enough reserves to meet all needs. In this case, the problem can be reduced to the transport by way of "cutting" applications.

Similar to the previous case, you can enter the value of $A_F$ (fictitious point of departure) and represent the difference:

$$A_f = \sum_{i=1}^{m} a_i - \sum_{j=1}^{n} b_j$$

## 6.4. Estimation of production volume using the Cobb-Douglas production function

We introduce such a thing as "production function". This is a function that displays the relationship between the factors of production and the maximum possible volume of the product. Also using this function, you can determine the minimum amount of costs required to produce a product for a given volume.

The main properties that this function has:

1. Due to the increase in costs for one resource, the production volume may increase to a specific value (it is not possible to hire workers in quantity more than in places in the room).

2. The complementarity and interchangeability of production factors.

In formal form, the production function can be represented as:

$$Q = f(K, L, M, T, N),$$

where $K$ – capital, $L$ – labor, $M$ – raw materials, $T$ – technologies, $N$ – entrepreneurial abilities.

We consider the two-factor model of the Cobb-Douglas production function because of its simplicity. With this model, you can reveal the relationship of labor ($T$) and capital ($K$). You can see that these factors are complementary and interchangeable. The model is formed as follows:

$$Q = AK^{\alpha}L^{\beta},$$

where $A$ is the coefficient of production (reflects the proportionality of functions and changes as the base technology changes), $K$ is capital, $L$ is labor, $\alpha, \beta$ are the elasticities of production for capital and labor, respectively.

In the production function $Q$ can be distinguished:

1. If $\alpha + \beta = 1$, then the function is proportionally increasing.

2. If $\alpha + \beta > 1$, the function is disproportionately increasing.

3. If $\alpha + \beta < 1$, the function decreases.

Let some firm have a short period of time of activity. Of the two factors for the variable we take the labor. This means that the firm can increase production by using more labor resources. Below is a graph of the Cobb - Douglas production function (Fig. 1). On it you can select the curve of the production function with one variable - $TP_L$

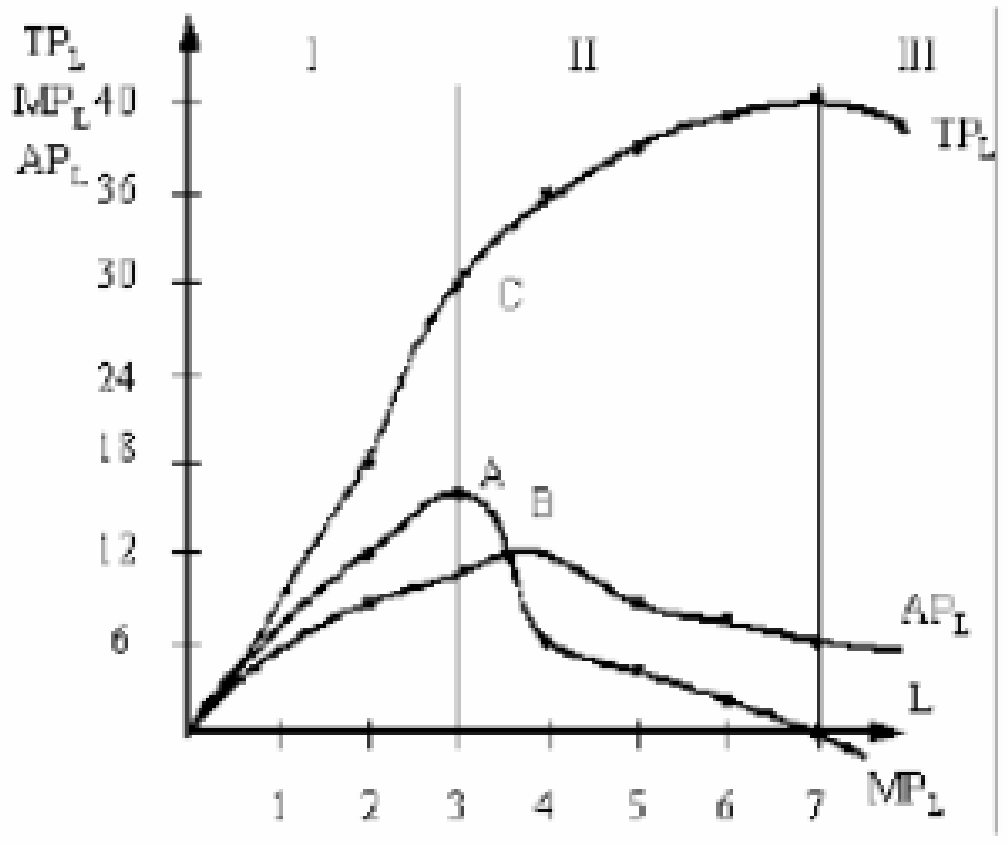

Fig 1. Dynamics and interrelation of total average and marginal products

We fix that in the short term there is a law of diminishing marginal productivity. An important condition is that one production factor is unchanged for a given time. Technique and production technology also remain unchanged. If the latest inventions are applied, the increase in output will be achieved using the same production factors.

Consider the case when capital is a fixed factor, and labor is variable. In this case, the firm increases production by using more labor resources. If you follow the law of diminishing marginal productivity, it is clear that a consistent increase in the variable resource leads to a diminishing return of this factor. It means that there is a decrease in the marginal product or the marginal productivity of labor. In the case of continuing the hiring of workers, the marginal productivity will take a negative value, which means that the volume of output will decrease.

We introduce the concept of "marginal productivity of labor" (in other words, "marginal product of labor"):

$$MP_L = \frac{\Delta Q_L}{\Delta L}$$

where $\Delta Q_L$ is the increase in production, $\Delta L$ is the increase in labor. Otherwise, this formula can be represented as:

$$MP_L = \frac{\Delta TP_L}{\Delta L}$$

where $\Delta TP_L$ is the performance gain to the total product.

$MP_K$ is defined in the same way. According to the schedule, you can analyze the relationship of total ($TP_L$), average ($AP_L$) and marginal products ($MP_L$).

There are three main stages in the movement of the total product curve:

1. The product limit increases (each new worker brings more production), therefore *TP* reaches a maximum at point *A*. Here the growth rate of the function is maximum.

2. The *MP* curve falls due to the law of diminishing returns. It can be seen that the growth rate of *TP* after *TC* slows down. $MP > 0$. When $MP = 0$, *TP* reaches maximum

3. *MP* is negative, and therefore *TP* begins to decline. This is due to the fact that the number of workers becomes redundant in relation to fixed capital.

The dynamics of the *MP* curve determine the configuration of the average product curve *AP*. It is seen that at the first stage, these curves increase, while the increment of output from newly hired workers is more than the average productivity of previously hired workers ($AP_L$). After point *A*, the average sample of four workers is reduced, as the fourth worker adds to the total product *TP* less than the third worker. There is a "scale effect":

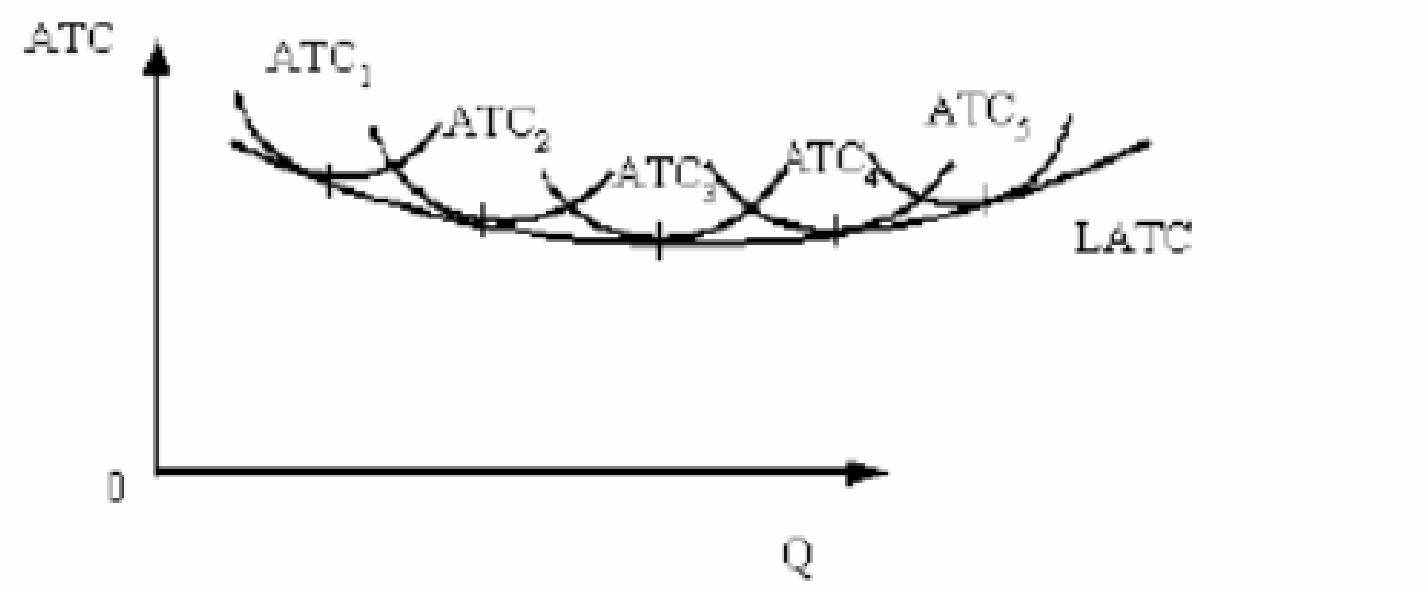

Fig 2. The curve of long-term and average costs of the company

1. It is appear in the lasting average cost of production (LATC)

2. LATC curve - envelope of the minimum short-term average cost of the company (per unit of output)

3. Changes in the number of all production factors in the long-term period in the activities of the company

When the company scale are changing, the values LATC may vary (Fig. 3)

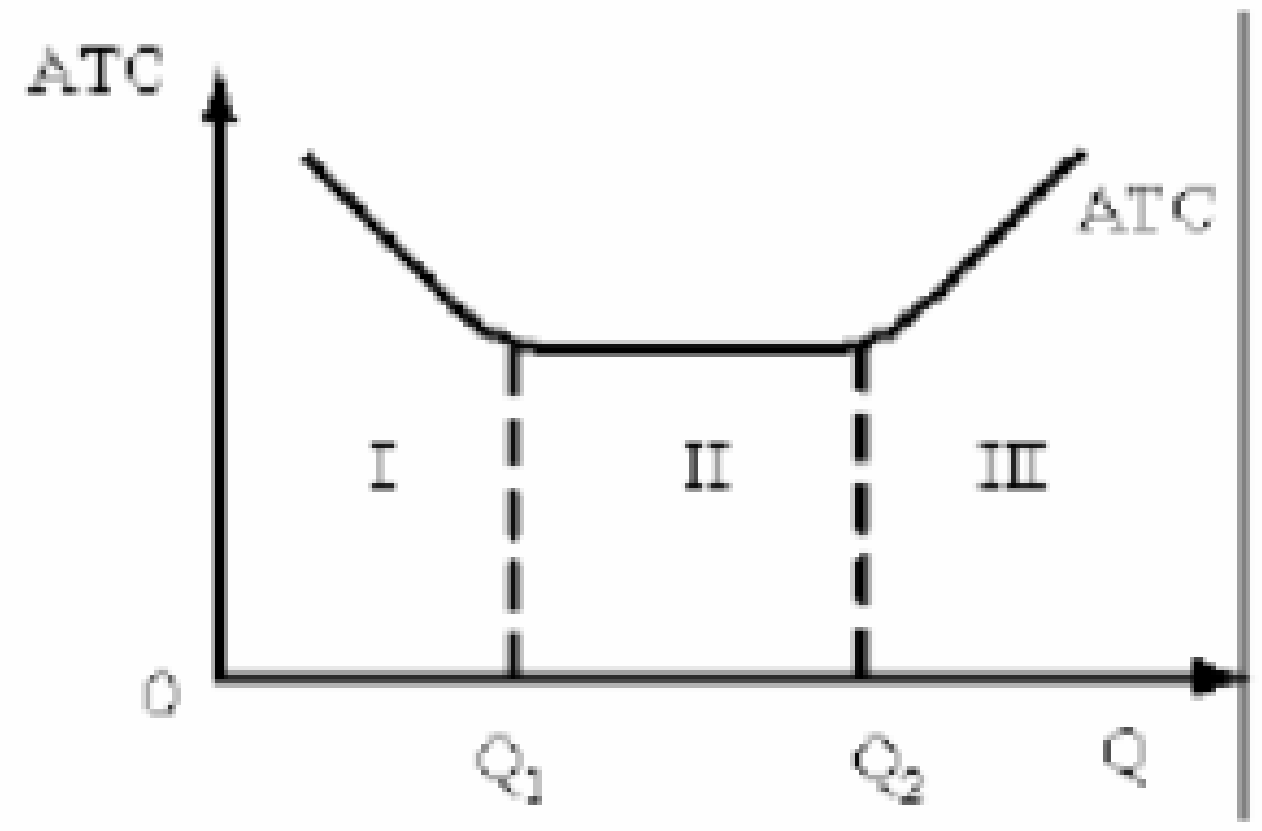

Fig 3.Dynamics of long-term average costs

There are the main steps:

**I. The positive effect of scale**

An increase in output decreases LATC. This is due to the saving effect.

**II. Constant returns to scale**

Costs remain constant when amount changes. The increase in the number of used resources causes an increase in production volumes.

**III. Negative scale effect**

The increase in production entails an increase in LATC. The cause of damage can be either technical factors (large enterprise size) or any organizational reasons (for example, the growth and inflexibility of the administrative and management apparatus).

## 8. Example

We would like to consider an example of the construction of the described model of multi-agent interaction in the mathematical problem of territorial distribution of 2 points of raw materials extraction, 2 points of temporary storage of raw materials, 2 points of production, 2 points of temporary storage of manufactured products, 4 stores.

In the problem at hand, a network is defined on the plane containing 18 vertices with coordinates $x_1,..,x_{18}$.

In the network nodes with coordinates $x_1$ and $x_6$ there are two points of extraction of raw materials $A_1$ and $A_2$, each of which extracts two types of raw materials $(\alpha_1, \alpha_2)$. The cost of extracting a unit of raw materials in warehouses can be expressed in the following table:

| | $\alpha_1$ | $\alpha_2$ |
|---|---|---|
| $A_1(x_1)$ | 1 | |
| $A_2(x_6)$ | | 2 |

The network can have two points of temporary storage of raw materials $S_1'$ and $S_2'$, and in two of the four points ($x_2, x_3, x_4, x_5$). 15 units are taken for storage services for the $\alpha_1$ type raw material and 22 units for the $\alpha_2$ type raw material.

The network can have two points of production of the final product B, and two of the four points ($x_7, x_{12}, x_{13}, x_{18}$) set the costs of raw materials for the production of each unit of manufactured products, which can be displayed in the following table:

| | $\alpha_1$ | $\alpha_2$ |
|---|---|---|
| $\beta_1$ | 1 | 1 |
| $B_2$ | 1 | 2 |
| $\beta_3$ | 2 | 1 |

Two points of temporary storage of final products $S_1''$ and $S_2''$ can be located on the network, and in two out of four points ($x_8$, $x_9$, $x_{10}$, $x_{11}$). For storage services of manufactured products, he takes 19 units for products of the type $\beta_1$, 33 units for products of the type $\beta_2$, 22 units for products of the type $\beta_3$.

The network has four stores $M_1$, $M_2$, $M_3$ and $M_4$, in which the manufactured products are sold at points ($x_{14}$, $x_{15}$, $x_{16}$, $x_{17}$), and at each point the demand for final products is determined, which can be displayed as follows:

| | $\beta_1$ | $\beta_2$ | $\beta_3$ |
|---|---|---|---|
| $X_{14}$ | 5 | 4 | 3 |
| $X_{15}$ | 5 | 3 | 4 |
| $X_{16}$ | 4 | 5 | 4 |
| $X_{17}$ | 3 | 5 | 5 |

Then the total demand:

| $\beta_1$ | $\beta_2$ | $\beta_3$ |
|---|---|---|
| 17 | 17 | 16 |

The total need for raw materials:

| | $\alpha_1$ | $\alpha_2$ |
|---|---|---|
| | 17*1+17*1+16*2 | 17*1+17*2+16*1 |
| Total | 66 | 67 |

Let us suppose the costs of transporting a unit of raw materials and a unit of manufactured products in horizontal and vertical directions will be represented as follows:

| | $\alpha_1$ | $\alpha_2$ | $\beta_1$ | $\beta_2$ | $\beta_3$ |
|---|---|---|---|---|---|
| Horizontally | 1 | 2 | 1 | 2 | 2 |
| Vertically | 2 | 1 | 1 | 2 | 1 |

The total costs per unit of manufactured products are equal to the sum of the costs of extracting the necessary raw materials, transporting them to the point of temporary storage of raw materials, storage costs for storing raw materials, transporting raw materials to the point of production of manufactured products, producing manufactured products, transporting manufactured productsto the point of temporary storage of manufactured products, storage costs for the storage of manufactured products, transportation to stores in which the sale of final products.

There are three agents who own different production sites. $Agent_3$ owns stores, $agent_2$ owns interchanges and warehouses, $agent_1$owns all production facilities. Our issue is to find a compromise solution in the problem, when everyone tries to maximize his net profit. $Agent_1$ costs: transport and storage costs 20% of the value of the stored product. Thus, its net profit is equal to the difference in income from storage and the costs of transporting and storing raw materials and manufactured products. Total storage income:

| | $\alpha_1$ | $\alpha_2$ | $\beta_1$ | $\beta_2$ | $\beta_3$ |
|---|---|---|---|---|---|
| 1 unit | 15,00 | 22,00 | 19,00 | 33,00 | 22,00 |
| The required number of units | 66,00 | 67,00 | 17,00 | 17,00 | 16,00 |
| Income from storing all units | 990,00 | 1474,00 | 323,00 | 561,00 | 352,00 |
| Total: | 2464,00 | | 1236,00 | | |
| The overall result: | 3700,00 | | | | |

One raw material warehouse can serve exactly one factory. Calculate the costs of all types of transportation. For this we use algorithm of Floyd and calculate the weight of the shortest paths in the graph. For transportation from the points of extraction of raw materials to the places of the proposed points of storage of raw materials:

| | $A_1$ ($x_1$) | $A_2$ ($x_6$) |
|---|---|---|
| | $\alpha_1$ | $\alpha_2$ |
| $X_2$ | 1 | 8 |
| $X_3$ | 2 | 6 |
| $X_4$ | 3 | 4 |
| $X_5$ | 4 | 2 |

For transportation from the locations of the proposed points of storage of raw materials to the places of intended production:

| | $S'(x_2)$ | | $S'(x_3)$ | | $S'(x_4)$ | | $S'(x_5)$ | |
|---|---|---|---|---|---|---|---|---|
| | $\alpha_1$ | $\alpha_2$ | $\alpha_1$ | $\alpha_2$ | $\alpha_1$ | $\alpha_2$ | $\alpha_1$ | $\alpha_2$ |
| $X_7$ | 3 | 3 | 4 | 5 | 5 | 7 | 6 | 9 |
| $X_{12}$ | 6 | 9 | 5 | 7 | 4 | 5 | 3 | 3 |
| $X_{13}$ | 5 | 4 | 6 | 6 | 7 | 8 | 8 | 10 |
| $X_{18}$ | 7 | 10 | 6 | 8 | 5 | 6 | 4 | 4 |

For transportation from the places of intended production to the places of proposed storage points for manufactured products:

| | $S''(x_{14})$ | | | $S''(x_{15})$ | | | $S''(x_{16})$ | | | $S''(x_{17})$ | | |
|---|---|---|---|---|---|---|---|---|---|---|---|---|
| | $\beta_1$ | $\beta_2$ | $\beta_3$ | $\beta_1$ | $\beta_2$ | $\beta_3$ | $\beta_1$ | $\beta_2$ | $\beta_3$ | $\beta_1$ | $\beta_2$ | $\beta_3$ |
| $X_7$ | 1 | 2 | 2 | 4 | 8 | 8 | 2 | 4 | 3 | 5 | 10 | 9 |
| $X_{12}$ | 2 | 4 | 4 | 3 | 6 | 6 | 3 | 6 | 5 | 4 | 8 | 7 |
| $X_{13}$ | 3 | 6 | 6 | 2 | 4 | 4 | 4 | 8 | 7 | 3 | 6 | 5 |
| $X_{18}$ | 4 | 8 | 8 | 1 | 2 | 2 | 5 | 10 | 9 | 2 | 4 | 3 |

For transportation from the places of intended production to the places of proposed storage points for manufactured products:

| | $M(x_8)$ | | | $M(x_9)$ | | | $M(x_{10})$ | | | $M(x_{11})$ | | |
|---|---|---|---|---|---|---|---|---|---|---|---|---|
| | $\beta_1$ | $\beta_2$ | $\beta_3$ | $\beta_1$ | $\beta_2$ | $\beta_3$ | $\beta_1$ | $\beta_2$ | $\beta_3$ | $\beta_1$ | $\beta_2$ | $\beta_3$ |
| $X_{14}$ | 1 | 2 | 1 | 2 | 4 | 3 | 3 | 6 | 5 | 4 | 8 | 7 |
| $X_{15}$ | 2 | 4 | 3 | 1 | 2 | 1 | 2 | 4 | 3 | 3 | 6 | 5 |
| $X_{16}$ | 3 | 6 | 5 | 2 | 4 | 3 | 1 | 2 | 1 | 2 | 4 | 3 |
| $X_{17}$ | 4 | 8 | 7 | 3 | 6 | 5 | 2 | 4 | 3 | 1 | 2 | 1 |

For transportation from the points of extraction of raw materials to the places of proposed points of production of manufactured products:

| | $x_1$ | | | | | | | | | | | | | | | |
|---|---|---|---|---|---|---|---|---|---|---|---|---|---|---|---|---|
| | $X_2$ | | | | $X_3$ | | | | $X_4$ | | | | $X_5$ | | | |
| | $X_7$ | $X_{12}$ | $X_{13}$ | $X_{18}$ | $X_7$ | $X_{12}$ | $X_{13}$ | $X_{18}$ | $X_7$ | $X_{12}$ | $X_{13}$ | $X_{18}$ | $X_7$ | $X_{12}$ | $X_{13}$ | $X_{18}$ |
| $\alpha_1$ | 4 | 7 | 6 | 8 | 6 | 7 | 8 | 8 | 8 | 7 | 10 | 8 | 10 | 7 | 12 | 8 |
| | $X_6$ | | | | | | | | | | | | | | | |
| | $X_2$ | | | | $X_3$ | | | | $X_4$ | | | | $X_5$ | | | |
| | $X_7$ | $X_{12}$ | $X_{13}$ | $X_{18}$ | $X_7$ | $X_{12}$ | $X_{13}$ | $X_{18}$ | $X_7$ | $X_{12}$ | $X_{13}$ | $X_{18}$ | $X_7$ | $X_{12}$ | $X_{13}$ | $X_{18}$ |
| $\alpha_2$ | 11 | 17 | 12 | 16 | 11 | 13 | 12 | 14 | 11 | 9 | 12 | 10 | 11 | 5 | 12 | 6 |

For transportation from the estimated points of production of manufactured products to stores where the sale of manufactured products:

| | $X_7$ | | | | | | | | | | | | | | | |
|---|---|---|---|---|---|---|---|---|---|---|---|---|---|---|---|---|
| | $X_8$ | | | | $X_9$ | | | | $X_{10}$ | | | | $X_{11}$ | | | |
| | $X_{14}$ | $X_{15}$ | $X_{16}$ | $X_{17}$ | $X_{14}$ | $X_{15}$ | $X_{16}$ | $X_{17}$ | $X_{14}$ | $X_{15}$ | $X_{16}$ | $X_{17}$ | $X_{14}$ | $X_{15}$ | $X_{16}$ | $X_{17}$ |
| $\beta_1$ | 2 | 3 | 4 | 5 | 4 | 3 | 4 | 5 | 6 | 5 | 4 | 5 | 8 | 7 | 6 | 5 |
| $\beta_2$ | 4 | 6 | 8 | 10 | 8 | 6 | 8 | 10 | 12 | 10 | 8 | 10 | 16 | 14 | 12 | 10 |
| $\beta_3$ | 3 | 5 | 7 | 9 | 7 | 5 | 7 | 9 | 11 | 9 | 7 | 9 | 15 | 13 | 11 | 9 |

| | $X_{12}$ | | | | | | | | | | | | | | | |
|---|---|---|---|---|---|---|---|---|---|---|---|---|---|---|---|---|
| | $X_8$ | | | | $X_9$ | | | | $X_{10}$ | | | | $X_{11}$ | | | |
| | $X_{14}$ | $X_{15}$ | $X_{16}$ | $X_{17}$ | $X_{14}$ | $X_{15}$ | $X_{16}$ | $X_{17}$ | $X_{14}$ | $X_{15}$ | $X_{16}$ | $X_{17}$ | $X_{14}$ | $X_{15}$ | $X_{16}$ | $X_{17}$ |
| $\beta_1$ | 5 | 6 | 7 | 8 | 5 | 4 | 5 | 6 | 5 | 4 | 3 | 4 | 5 | 4 | 3 | 2 |
| $\beta_2$ | 10 | 12 | 14 | 16 | 10 | 8 | 10 | 12 | 10 | 8 | 6 | 10 | 10 | 8 | 6 | 4 |
| $\beta_3$ | 9 | 11 | 13 | 15 | 9 | 7 | 9 | 11 | 9 | 7 | 5 | 7 | 9 | 7 | 5 | 3 |
| | $X_{13}$ | | | | | | | | | | | | | | | |
| | $X_8$ | | | | $X_9$ | | | | $X_{10}$ | | | | $X_{11}$ | | | |
| | $X_{14}$ | $X_{15}$ | $X_{16}$ | $X_{17}$ | $X_{14}$ | $X_{15}$ | $X_{16}$ | $X_{17}$ | $X_{14}$ | $X_{15}$ | $X_{16}$ | $X_{17}$ | $X_{14}$ | $X_{15}$ | $X_{16}$ | $X_{17}$ |
| $\beta_1$ | 3 | 4 | 5 | 6 | 5 | 4 | 5 | 6 | 7 | 6 | 5 | 6 | 9 | 8 | 7 | 6 |
| $\beta_2$ | 6 | 8 | 10 | 12 | 10 | 8 | 10 | 12 | 14 | 12 | 1 | 12 | 18 | 16 | 14 | 12 |
| $\beta_3$ | 4 | 6 | 8 | 10 | 8 | 6 | 8 | 10 | 12 | 10 | 8 | 10 | 16 | 14 | 12 | 10 |
| | $X_{18}$ | | | | | | | | | | | | | | | |
| | $X_8$ | | | | $X_9$ | | | | $X_{10}$ | | | | $X_{11}$ | | | |
| | $X_{14}$ | $X_{15}$ | $X_{16}$ | $X_{17}$ | $X_{14}$ | $X_{15}$ | $X_{16}$ | $X_{17}$ | $X_{14}$ | $X_{15}$ | $X_{16}$ | $X_{17}$ | $X_{14}$ | $X_{15}$ | $X_{16}$ | $X_{17}$ |
| $\beta_1$ | 6 | 7 | 8 | 9 | 6 | 5 | 6 | 7 | 6 | 5 | 4 | 4 | 6 | 5 | 4 | 3 |
| $\beta_2$ | 12 | 14 | 16 | 18 | 12 | 10 | 12 | 14 | 12 | 10 | 8 | 10 | 12 | 10 | 8 | 6 |
| $\beta_3$ | 10 | 12 | 14 | 16 | 10 | 8 | 10 | 12 | 10 | 8 | 6 | 8 | 10 | 8 | 6 | 4 |

For second agent, income is determined using the Cobb-Douglas production function, which is:

$$Q = AK^{\alpha}L^{\beta}$$

In our problem, we assume that the production function increases proportionally. Suppose that $K$, $L$ — costs for raw materials of the form $\alpha_1$ and $\alpha_2$, respectively, their use is expressed in power equivalent.

| | $\alpha_1$ | $\alpha_2$ |
|---|---|---|
| $\beta_1$ | 0,5 | 0,5 |
| $\beta_2$ | 0,33 | 0,33 |
| $\beta_3$ | 0,67 | 0,67 |

The production factor for each item of raw materials for each type of manufactured products can be expressed in the following table:

| | $X_7$ | $X_{12}$ | $X_{13}$ | $X_{18}$ |
|---|---|---|---|---|
| $\beta_1$ | 2,1 | 2,3 | 2,2 | 2,5 |
| $\beta_2$ | 2,2 | 2,2 | 2 | 2,1 |
| $\beta_3$ | 2,3 | 2,1 | 2,3 | 2,15 |

Each production point cannot produce more than 10 items of each type, which means that at the second we will produce what remains to meet the demand. Thus, the distribution of production in factories can be expressed in the following table (if the factory produces 10 units of goods of the type $\beta_1$, then the second does not necessarily have the same type of $\beta_2$, that is, the table is schematic):

| | First plant | Second plant |
|---|---|---|
| $\beta_1$ | 10 | 7 |
| $\beta_2$ | 10 | 7 |
| $\beta_3$ | 10 | 6 |

Then, we can calculate the production function at each possible location of the plant, and, accordingly, the net profit.

| $X_7$ | `Production factor | Amount $\alpha_1$ | Cost $\alpha_1$ | Total cost $\alpha_1$ | Power Equiva lent | Amount $\alpha_2$ | Cost $\alpha_2$ | Total $\alpha_2$ | Power Equiva lent | Total cost $\alpha_1$and$\alpha_2$ | Production function per unit of production | Production function for the entire release | Net profit |
|---|---|---|---|---|---|---|---|---|---|---|---|---|---|
| $\beta_1$ | 2,10 | 10,00 | 15,00 | 150,00 | 0,50 | 10,00 | 22,00 | 220,00 | 0,50 | 370,00 | 38,15 | 381,48 | 11,48 |
| $\beta_2$ | 2,20 | 10,00 | 15,00 | 150,00 | 0,33 | 20,00 | 22,00 | 440,00 | 0,67 | 590,00 | 67,87 | 678,65 | 88,65 |
| $\beta_3$ | 2,30 | 20,00 | 15,00 | 300,00 | 0,67 | 10,00 | 22,00 | 220,00 | 0,33 | 520,00 | 62,29 | 622,87 | 102,87 |

| $X_{12}$ | `Production factor | Amount $\alpha_1$ | Cost $\alpha_1$ | Total cost $\alpha_1$ | Power Equiva lent | Amount $\alpha_2$ | Cost $\alpha_2$ | Total $\alpha_2$ | Power Equiva lent | Total cost $\alpha_1$and$\alpha_2$ | Production function per unit of production | Production function for the entire release | Net profit |
|---|---|---|---|---|---|---|---|---|---|---|---|---|---|
| $\beta_1$ | 2,30 | 10,00 | 15,00 | 150,00 | 0,50 | 10,00 | 22,00 | 220,00 | 0,50 | 370,00 | 41,78 | 417,82 | 47,82 |
| $\beta_2$ | 2,20 | 10,00 | 15,00 | 150,00 | 0,33 | 20,00 | 22,00 | 440,00 | 0,67 | 590,00 | 67,87 | 678,65 | 88,65 |
| $\beta_3$ | 2,10 | 20,00 | 15,00 | 300,00 | 0,67 | 10,00 | 22,00 | 220,00 | 0,33 | 520,00 | 56,87 | 568,71 | 48,71 |

| $X_{13}$ | `Production factor | Amount $\alpha_1$ | Cost $\alpha_1$ | Total cost $\alpha_1$ | Power Equiva lent | Amount $\alpha_2$ | Cost $\alpha_2$ | Total $\alpha_2$ | Power Equiva lent | Total cost $\alpha_1$and$\alpha_2$ | Production function per unit of production | Production function for the entire release | Net profit |
|---|---|---|---|---|---|---|---|---|---|---|---|---|---|
| $\beta_1$ | 2,20 | 10,00 | 15,00 | 150,00 | 0,50 | 10,00 | 22,00 | 220,00 | 0,50 | 370,00 | 39,96 | 399,65 | 29,65 |
| $\beta_2$ | 2,00 | 10,00 | 15,00 | 150,00 | 0,33 | 20,00 | 22,00 | 440,00 | 0,67 | 590,00 | 61,70 | 616,95 | 26,95 |
| $\beta_3$ | 2,30 | 20,00 | 15,00 | 300,00 | 0,67 | 10,00 | 22,00 | 220,00 | 0,33 | 520,00 | 62,29 | 622,87 | 102,87 |

| $X_{18}$ | `Production factor | Amount $\alpha_1$ | Cost $\alpha_1$ | Total cost $\alpha_1$ | Power Equiva lent | Amount $\alpha_2$ | Cost $\alpha_2$ | Total $\alpha_2$ | Power Equiva lent | Total cost $\alpha_1$and$\alpha_2$ | Production function per unit of production | Production function for the entire release | Net profit |
|---|---|---|---|---|---|---|---|---|---|---|---|---|---|
| $\beta_1$ | 2,50 | 10,00 | 15,00 | 150,00 | 0,50 | 10,00 | 22,00 | 220,00 | 0,50 | 370,00 | 45,41 | 454,15 | 84,15 |
| $\beta_2$ | 2,10 | 10,00 | 15,00 | 150,00 | 0,33 | 20,00 | 22,00 | 440,00 | 0,67 | 590,00 | 64,78 | 647,80 | 57,80 |
| $\beta_3$ | 2,15 | 20,00 | 15,00 | 300,00 | 0,67 | 10,00 | 22,00 | 220,00 | 0,33 | 520,00 | 58,22 | 582,25 | 62,25 |

| $X_7$ | `Production factor | Amount $\alpha_1$ | Cost $\alpha_1$ | Total cost $\alpha_1$ | Power Equiva lent | Amount $\alpha_2$ | Cost $\alpha_2$ | Total $\alpha_2$ | Power Equiva lent | Total cost $\alpha_1$and$\alpha_2$ | Production function per unit of production | Production function for the entire release | Net profit |
|---|---|---|---|---|---|---|---|---|---|---|---|---|---|
| $\beta_1$ | 2,10 | 7,00 | 15,00 | 105,00 | 0,50 | 7,00 | 22,00 | 154,00 | 0,50 | 259,00 | 38,15 | 267,04 | 8,04 |
| $\beta_2$ | 2,20 | 7,00 | 15,00 | 105,00 | 0,33 | 14,00 | 22,00 | 308,00 | 0,67 | 413,00 | 67,87 | 475,06 | 62,06 |
| $\beta_3$ | 2,30 | 12,00 | 15,00 | 180,00 | 0,67 | 6,00 | 22,00 | 132,00 | 0,33 | 312,00 | 62,29 | 373,72 | 61,72 |

| $X_{12}$ | `Production factor | Amount $\alpha_1$ | Cost $\alpha_1$ | Total cost $\alpha_1$ | Power Equiva lent | Amount $\alpha_2$ | Cost $\alpha_2$ | Total $\alpha_2$ | Power Equiva lent | Total cost $\alpha_1$and$\alpha_2$ | Production function per unit of production | Production function for the entire release | Net profit |
|---|---|---|---|---|---|---|---|---|---|---|---|---|---|
| $\beta_1$ | 2,30 | 7,00 | 15,00 | 105,00 | 0,50 | 7,00 | 22,00 | 154,00 | 0,50 | 259,00 | 41,78 | 292,47 | 33,47 |
| $\beta_2$ | 2,20 | 7,00 | 15,00 | 105,00 | 0,33 | 14,00 | 22,00 | 308,00 | 0,67 | 413,00 | 67,87 | 475,06 | 62,06 |
| $\beta_3$ | 2,10 | 12,00 | 15,00 | 180,00 | 0,67 | 6,00 | 22,00 | 132,00 | 0,33 | 312,00 | 56,87 | 341,23 | 29,23 |

| $X_{13}$ | `Production factor | Amount $\alpha_1$ | Cost $\alpha_1$ | Total cost $\alpha_1$ | Power Equivalent | Amount $\alpha_2$ | Cost $\alpha_2$ | Total $\alpha_2$ | Power Equivalent | Total cost $\alpha_1$and$\alpha_2$ | Production function per unit of production | Production function for the entire release | Net profit |
|---|---|---|---|---|---|---|---|---|---|---|---|---|---|
| $\beta_1$ | 2,20 | 7,00 | 15,00 | 105,00 | 0,50 | 7,00 | 22,00 | 154,00 | 0,50 | 259,00 | 39,96 | 279,75 | 20,75 |
| $\beta_2$ | 2,00 | 7,00 | 15,00 | 105,00 | 0,33 | 14,00 | 22,00 | 308,00 | 0,67 | 413,00 | 61,70 | 431,87 | 18,87 |
| $\beta_3$ | 2,30 | 12,00 | 15,00 | 180,00 | 0,67 | 6,00 | 22,00 | 132,00 | 0,33 | 312,00 | 62,29 | 373,72 | 61,72 |

| $X_{18}$ | `Production factor | Amount $\alpha_1$ | Cost $\alpha_1$ | Total cost $\alpha_1$ | Power Equivalent | Amount $\alpha_2$ | Cost $\alpha_2$ | Total $\alpha_2$ | Power Equivalent | Total cost $\alpha_1$and$\alpha_2$ | Production function per unit of production | Production function for the entire release | Net profit |
|---|---|---|---|---|---|---|---|---|---|---|---|---|---|
| $\beta_1$ | 2,50 | 7,00 | 15,00 | 105,00 | 0,50 | 7,00 | 22,00 | 154,00 | 0,50 | 259,00 | 45,41 | 317,90 | 58,90 |
| $\beta_2$ | 2,10 | 7,00 | 15,00 | 105,00 | 0,33 | 14,00 | 22,00 | 308,00 | 0,67 | 413,00 | 64,78 | 453,46 | 40,46 |
| $\beta_3$ | 2,15 | 12,00 | 15,00 | 180,00 | 0,67 | 6,00 | 22,00 | 132,00 | 0,33 | 312,00 | 58,22 | 349,35 | 37,35 |

For the third agent, income is determined from the difference in profit from the sale of a unit of goods and its purchase price, which is determined by the production function at the plant, and the storage costs.

Prices in the store are determined, so we can calculate the total income from all points:

| | $\beta_1$ | $\beta_2$ | $\beta_3$ |
|---|---|---|---|
| Unit sales price | 85 | 125 | 115 |
| The required number of units | 17 | 17 | 16 |
| Revenue from the sale of all units | 1445 | 2125 | 1840 |
| Total | 5410 | | |

Thus, we deal with six situations on which all actions of players depend - depending on the location of production, we have temporary storage facilities and shops.

1.($X_7$,$X_{12}$)

Production is distributed as follows.

| **Production revenue** | **At the $X_7$ factory** | **Net profit at the factory $X_7$** | **At the $X_{12}$ factory** | **Net profit at the factory $X_{12}$** | **Total:** |
|---|---|---|---|---|---|
| $\beta_1$ | 7 | 8,04 | 10 | 47,82 | 55,85 |
| $\beta_2$ | 10 | 88,65 | 7 | 62,06 | 150,71 |
| $\beta_3$ | 10 | 102,87 | 6 | 29,23 | 132,10 |
| Total: | 27 | 199,56 | 23 | 139,11 | 338,66 |

Net profit of stores will be:

| | The cost of the factory $X_7$ | Quantity of goods | Storage costs | Total Cost Per Unit | Total | The cost of the factory $X_{12}$ | Quantity of goods | Storage costs | Total Cost Per Unit | Total | **Net profit** |
|---|---|---|---|---|---|---|---|---|---|---|---|
| $\beta_1$ | 38,15 | 7,00 | 19,00 | 57,15 | 400,04 | 41,78 | 10,00 | 19,00 | 60,78 | 607,82 | **437,15** |
| $\beta_2$ | 67,87 | 10,00 | 33,00 | 100,87 | 1008,65 | 67,87 | 7,00 | 33,00 | 100,87 | 706,06 | **410,29** |
| $\beta_3$ | 62,29 | 10,00 | 22,00 | 84,29 | 842,87 | 56,87 | 6,00 | 22,00 | 78,87 | 473,23 | **523,90** |
| Total | | | | | | | | | | 4038,66 | **1371,34** |

Consider the required amount of each manufacture resources:

| Factory $x_7$ | | | | | | Factory $x_{12}$ | | | | | |
|---|---|---|---|---|---|---|---|---|---|---|---|
| $\beta_1$ | | $\beta_2$ | | $\beta_3$ | | $\beta_1$ | | $\beta_2$ | | $\beta_3$ | |
| 7 | | 10 | | 10 | | 10 | | 7 | | 6 | |
| The required amount of resource per unit of production | | | | | | | | | | | |
| $\alpha_1$ | $\alpha_2$ | $\alpha_1$ | $\alpha_2$ | $\alpha_1$ | $\alpha_2$ | $\alpha_1$ | $\alpha_2$ | $\alpha_1$ | $\alpha_2$ | $\alpha_1$ | $\alpha_2$ |
| 1 | 1 | 1 | 2 | 2 | 1 | 1 | 1 | 1 | 2 | 2 | 1 |
| The required amount of resource for the entire quantity of products | | | | | | | | | | | |
| 7 | 7 | 10 | 20 | 20 | 10 | 10 | 10 | 7 | 14 | 12 | 6 |
| Total resources $\alpha_1$ | | | Total resources $\alpha_2$ | | | Total resources $\alpha_1$ | | | Total resources $\alpha_2$ | | |
| 37 | | | 37 | | | 29 | | | 30 | | |

The cost of transporting a unit of production from the points of extraction of raw materials to the points of production:

| | $x_2$ | | $x_3$ | | $x_4$ | | $x_5$ | |
|---|---|---|---|---|---|---|---|---|
| | $x_7$ | $x_{12}$ | $x_7$ | $x_{12}$ | $x_7$ | $x_{12}$ | $x_7$ | $x_{12}$ |
| $\alpha_1$ | 4 | 7 | 6 | 7 | 8 | 7 | 10 | 7 |
| $\alpha_2$ | 11 | 17 | 11 | 13 | 11 | 9 | 11 | 5 |

It is seen that it is most advantageous to transport through points $x_2$ and $x_5$. Then the profit minus the costs of transportation and storage of raw materials will be:

| | Factory $x_7$ | | | | Factory $x_{12}$ | | | |
|---|---|---|---|---|---|---|---|---|
| | **cost per unit of transportation** | **storage costs per unit** | **total units** | **total** | **cost per unit of transportation** | **storage costs per unit** | **total units** | **total** |
| $\alpha_1$ | 4,00 | 0,20 | 37,00 | 155,40 | 7,00 | 0,20 | 29,00 | 208,80 |
| $\alpha_2$ | 11,00 | 0,20 | 37,00 | 421,80 | 5,00 | 0,20 | 30,00 | 162,00 |
| Net profit | | | | | | | | 1516,00 |

The cost of transporting a unit of manufactured products from production points to stores where the sale of final products is carried out

| | $X_7$ | | | | | | | | | | | | | | | |
|---|---|---|---|---|---|---|---|---|---|---|---|---|---|---|---|---|
| | $X_8$ | | | | $X_9$ | | | | $X_{10}$ | | | | $X_{11}$ | | | |
| | $X_{14}$ | $X_{15}$ | $X_{16}$ | $X_{17}$ | $X_{14}$ | $X_{15}$ | $X_{16}$ | $X_{17}$ | $X_{14}$ | $X_{15}$ | $X_{16}$ | $X_{17}$ | $X_{14}$ | $X_{15}$ | $X_{16}$ | $X_{17}$ |
| $\beta_1$ | 2 | 3 | 4 | 5 | 4 | 3 | 4 | 5 | 6 | 5 | 4 | 5 | 8 | 7 | 6 | 5 |
| $\beta_2$ | 4 | 6 | 8 | 10 | 8 | 6 | 8 | 10 | 12 | 10 | 8 | 10 | 16 | 14 | 12 | 10 |
| $\beta_3$ | 3 | 5 | 7 | 9 | 7 | 5 | 7 | 9 | 11 | 9 | 7 | 9 | 15 | 13 | 11 | 9 |
| | $X_{12}$ | | | | | | | | | | | | | | | |
| | $X_8$ | | | | $X_9$ | | | | $X_{10}$ | | | | $X_{11}$ | | | |
| | $X_{14}$ | $X_{15}$ | $X_{16}$ | $X_{17}$ | $X_{14}$ | $X_{15}$ | $X_{16}$ | $X_{17}$ | $X_{14}$ | $X_{15}$ | $X_{16}$ | $X_{17}$ | $X_{14}$ | $X_{15}$ | $X_{16}$ | $X_{17}$ |
| $\beta_1$ | 5 | 6 | 7 | 8 | 5 | 4 | 5 | 6 | 5 | 4 | 3 | 4 | 5 | 4 | 3 | 2 |
| $\beta_2$ | 10 | 12 | 14 | 16 | 10 | 8 | 10 | 12 | 10 | 8 | 6 | 10 | 10 | 8 | 6 | 4 |
| $\beta_3$ | 9 | 11 | 13 | 15 | 9 | 7 | 9 | 11 | 9 | 7 | 5 | 7 | 9 | 7 | 5 | 3 |

Transportation is carried out according to the principle: "first, the cheapest transportation, after more expensive", (through warehouses $x_8$ and $x_{11}$ - according to a preliminary assessment of the previous table):

| | **Required amount $\beta_1$** | **From the factory $X_7$** | **From the factory $X_{12}$** | **Required amount $\beta_2$** | **From the factory $X_7$** | **From the factory $X_{12}$** | **Required amount $\beta_3$** | **From the factory $X_7$** | **From the factory $X_{12}$** |
|---|---|---|---|---|---|---|---|---|---|
| $X_{14}$ | 5 | 5 | 0 | 4 | 4 | 0 | 3 | 3 | 0 |
| $X_{15}$ | 5 | 2 | 3 | 3 | 3 | 0 | 4 | 4 | 0 |
| $X_{16}$ | 4 | 0 | 4 | 5 | 3 | 2 | 4 | 0 | 4 |
| $X_{17}$ | 3 | 0 | 3 | 5 | 4 | 1 | 5 | 0 | 5 |
| **Total shipping costs** | | | | | | | | | |
| $X_{14}$ | | 10 | 0 | | 16 | 0 | | 9 | 0 |
| $X_{15}$ | | 6 | 12 | | 18 | 0 | | 20 | 0 |
| $X_{16}$ | | 0 | 16 | | 24 | 12 | | 0 | 20 |
| $X_{17}$ | | 0 | 6 | | 40 | 4 | | 0 | 15 |
| **Total** | | | | | | | | | 228 |

Income from storage in warehouses:

| | **Cost of purchase at the factory $X_7$** | **Amount** | **Total** | **Cost of purchase at the factory $X_{12}$** | **Amount** | **Total** |
|---|---|---|---|---|---|---|
| $\beta_1$ | 38,15 | 7,00 | 267,04 | 41,78 | 10,00 | 417,82 |
| $\beta_2$ | 67,87 | 10,00 | 678,65 | 67,87 | 7,00 | 475,06 |
| $\beta_3$ | 62,29 | 10,00 | 622,87 | 56,87 | 6,00 | 341,23 |
| Total | | | | | | 2802,66 |
| Net profit | | | | | | 675,47 |

Thus, the net profit of the first agent will be: 1516 + 675.47-228 = **1963.47**

2)($X_7$,$X_{13}$)

Production will be distributed as follows:

| **Production revenue** | **At the $X_7$factory** | **Net profit at the factory $X_7$** | **At the $X_{13}$factory** | **Net profit at the factory $X_{13}$** | **Total:** |
|---|---|---|---|---|---|
| $\beta_1$ | 7 | 8,04 | 10 | 29,65 | 37,69 |
| $\beta_2$ | 10 | 88,65 | 7 | 18,87 | 107,52 |
| $\beta_3$ | 10 | 102,87 | 6 | 61,72 | 164,59 |
| Total: | 27 | 199,56 | 23 | 110.24 | 309,80 |

Net profit of stores will be:

| | The cost of the factory $X_7$ | Quantity of goods | Storage costs | Total Cost Per Unit | Total | The cost of the factory $X_{13}$ | Quantity of goods | Storage costs | Total Cost Per Unit | Total | **Net profit** |
|---|---|---|---|---|---|---|---|---|---|---|---|
| $\beta_1$ | 38,15 | 7,00 | 19,00 | 57,15 | 400,04 | 39,96 | 10,00 | 19,00 | 58,96 | 589,65 | 455,31 |
| $\beta_2$ | 67,87 | 10,00 | 33,00 | 100,87 | 1008,65 | 61,70 | 7,00 | 33,00 | 94,70 | 662,87 | 453,48 |
| $\beta_3$ | 62,29 | 10,00 | 22,00 | 84,29 | 842,87 | 62,29 | 6,00 | 22,00 | 84,29 | 505,72 | 491,41 |
| Total | | | | | | | | | | 4009,80 | 1400,20 |

Consider the required amount of resources for each production:

| Factory $x_7$ | | | | | | Factory $x_{13}$ | | | | | |
|---|---|---|---|---|---|---|---|---|---|---|---|
| $\beta_1$ | | $\beta_2$ | | $\beta_3$ | | $\beta_1$ | | $\beta_2$ | | $\beta_3$ | |
| 7 | | 10 | | 10 | | 10 | | 7 | | 6 | |
| The required amount of resource per unit of production | | | | | | | | | | | |
| $\alpha_1$ | $\alpha_2$ | $\alpha_1$ | $\alpha_2$ | $\alpha_1$ | $\alpha_2$ | $\alpha_1$ | $\alpha_2$ | $\alpha_1$ | $\alpha_2$ | $\alpha_1$ | $\alpha_2$ |
| 7 | 7 | 10 | 20 | 20 | 10 | 10 | 10 | 7 | 14 | 12 | 6 |
| The required amount of resource for the entire quantity of products | | | | | | | | | | | |
| 7 | 7 | 10 | 20 | 20 | 10 | 10 | 10 | 7 | 14 | 12 | 6 |
| Total resources $\alpha_1$ | | | Total resources $\alpha_2$ | | | Total resources $\alpha_1$ | | | Total resources $\alpha_2$ | | |
| 37 | | | 37 | | | 29 | | | 30 | | |

The cost of transporting a unit of production from the points of extraction of raw materials to the points of production:

| | $x_2$ | | $x_3$ | | $x_4$ | | $x_5$ | |
|---|---|---|---|---|---|---|---|---|
| | $x_7$ | $x_{13}$ | $x_7$ | $x_{13}$ | $x_7$ | $x_{13}$ | $x_7$ | $x_{13}$ |
| $\alpha_1$ | 4 | 6 | 6 | 8 | 8 | 10 | 10 | 12 |
| $\alpha_2$ | 11 | 12 | 11 | 12 | 11 | 12 | 11 | 12 |

It is seen that it is most advantageous to transport through points $x_2$ and $x_3$. Then the profit minus the costs of transportation and storage of raw materials will be:

| | Factory $x_7$ | | | | Factory $x_{13}$ | | | |
|---|---|---|---|---|---|---|---|---|
| | **cost per unit of transportation** | **storage costs per unit** | **total units** | **total** | **cost per unit of transportation** | **storage costs per unit** | **total units** | **total** |
| $\alpha_1$ | 4,00 | 0,20 | 37,00 | 162,80 | 8,00 | 0,20 | 29,00 | 243,60 |
| $\alpha_2$ | 11,00 | 0,20 | 37,00 | 421,80 | 12,00 | 0,20 | 30,00 | 372,00 |
| Net profit | | | | | | | | 1263,80 |

The cost of transporting a unit of manufactured products from production points to stores where the sale of final products is carried out, is carried out according to the principle: "first, the cheapest transportation, after more expensive", through warehouses $x_8$ and $x_{10}$:

| | **Required amount β₁** | **From the factory X₇** | **From the factory X₁₃** | **Required amount β₂** | **From the factory X₇** | **From the factory X₁₃** | **Required amount β₃** | **From the factory X₇** | **From thef actory X₁₃** |
|---|---|---|---|---|---|---|---|---|---|
| $X_{14}$ | 5 | 5 | 0 | 4 | 4 | 0 | 3 | 3 | 0 |
| $X_{15}$ | 5 | 2 | 3 | 3 | 3 | 0 | 4 | 4 | 0 |
| $X_{16}$ | 4 | 0 | 4 | 5 | 3 | 2 | 4 | 0 | 4 |
| $X_{17}$ | 3 | 0 | 3 | 5 | 4 | 1 | 5 | 0 | 5 |
| **Total shipping costs** | | | | | | | | | |
| $X_{14}$ | | 10 | 0 | | 16 | 0 | | 9 | 0 |
| $X_{15}$ | | 6 | 18 | | 18 | 0 | | 20 | 0 |
| $X_{16}$ | | 0 | 16 | | 24 | 2 | | 0 | 32 |
| $X_{17}$ | | 0 | 18 | | 40 | 12 | | 0 | 50 |
| **Total** | | | | | | | | | 291 |

Income from storage in warehouses:

| | **Cost of purchase at the factory X₇** | **Amount** | **Total** | **Cost of purchase at the factory X₁₂** | **Amount** | **Total** |
|---|---|---|---|---|---|---|
| $\beta_1$ | 38,15 | 7,00 | 267,04 | 39,96 | 10,00 | 399,65 |
| $\beta_2$ | 67,87 | 10,00 | 678,65 | 61,70 | 7,00 | 431,87 |
| $\beta_3$ | 62,29 | 10,00 | 622,87 | 62,29 | 6,00 | 373,72 |
| Total | | | | | | 2773,80 |
| Net profit | | | | | | 681,24 |

The net profit of the first agent will be **1654.04**

3) ($X_7$,$X_{18}$)

Production will be distributed as follows:

| Production revenue | At the $X_7$factory | Net profit at the factory $X_7$ | At the $X_{18}$factory | Net profit at the factory $X_{18}$ | Total: |
|---|---|---|---|---|---|
| $\beta_1$ | 7 | 8,04 | 10 | 84,15 | 92,19 |
| $\beta_2$ | 10 | 88,65 | 7 | 40,46 | 129,11 |
| $\beta_3$ | 10 | 102,87 | 6 | 37,35 | 140,22 |
| Total: | 27 | 199.56 | 23 | 161.96 | 361,52 |

Net profit of stores:

| | The cost of the factory $X_7$ | Quantity of goods | Storage costs | Total Cost Per Unit | Total | The cost of the factory $X_{18}$ | Quantity of goods | Storage costs | Total Cost Per Unit | Total | **Net profit** |
|---|---|---|---|---|---|---|---|---|---|---|---|
| $\beta_1$ | 38,15 | 7,00 | 19,00 | 57,15 | 400,04 | 45,41 | 10,00 | 19,00 | 64,41 | 644,15 | 400,81 |
| $\beta_2$ | 67,87 | 10,00 | 33,00 | 100,87 | 1008,65 | 64,78 | 7,00 | 33,00 | 97,78 | 684,46 | 431,89 |
| $\beta_3$ | 62,29 | 10,00 | 22,00 | 84,29 | 842,87 | 58,22 | 6,00 | 22,00 | 80,22 | 481,35 | 515,78 |
| Total | | | | | | | | | | 4061,52 | 1348,48 |

Consider the required amount of resources for each production:

| Factory $x_7$ | | | | | | Factory $x_{18}$ | | | | | |
|---|---|---|---|---|---|---|---|---|---|---|---|
| $\beta_1$ | | $\beta_2$ | | $\beta_3$ | | $\beta_1$ | | $\beta_2$ | | $\beta_3$ | |
| 7 | | 10 | | 10 | | 10 | | 7 | | 6 | |
| The required amount of resource per unit of production | | | | | | | | | | | |
| $\alpha_1$ | $\alpha_2$ | $\alpha_1$ | $\alpha_2$ | $\alpha_1$ | $\alpha_2$ | $\alpha_1$ | $\alpha_2$ | $\alpha_1$ | $\alpha_2$ | $\alpha_1$ | $\alpha_2$ |
| 1 | 1 | 1 | 2 | 2 | 1 | 1 | 1 | 1 | 2 | 2 | 1 |
| The required amount of resource for the entire quantity of products | | | | | | | | | | | |
| 7 | 7 | 10 | 20 | 20 | 10 | 10 | 10 | 7 | 14 | 12 | 6 |

| Total resources $\alpha_1$ | Total resources $\alpha_2$ | Total resources $\alpha_1$ | Total resources $\alpha_2$ |
|---|---|---|---|
| 37 | 37 | 29 | 30 |

The cost of transporting a unit of production from the points of extraction of raw materials to the points of production:

| | $x_2$ | | $x_3$ | | $x_4$ | | $x_5$ | |
|---|---|---|---|---|---|---|---|---|
| | $x_7$ | $x_{18}$ | $x_7$ | $x_{18}$ | $x_7$ | $x_{18}$ | $x_7$ | $x_{18}$ |
| $\alpha_1$ | 4 | 8 | 6 | 8 | 8 | 8 | 10 | 8 |
| $\alpha_2$ | 11 | 16 | 11 | 14 | 11 | 10 | 11 | 6 |

It is seen that it is most advantageous to transport through points $x_2$ and $x_5$. Then the profit minus the costs of transportation and storage of raw materials will be:

| | Factory $x_7$ | | | | Factory $x_{18}$ | | | |
|---|---|---|---|---|---|---|---|---|
| | **cost per unit of transportation** | **storage costs per unit** | **total units** | **total** | **cost per unit of transportation** | **storage costs per unit** | **total units** | **total** |
| $\alpha_1$ | 4,00 | 0,20 | 37,00 | 162,80 | 8,00 | 0,20 | 29,00 | 243,60 |
| $\alpha_2$ | 11,00 | 0,20 | 37,00 | 421,80 | 6,00 | 0,20 | 30,00 | 192,00 |
| Net profit | | | | | | | | 1443,80 |

The cost of transporting a unit of manufactured products from production points to stores where the sale of final products is carried out, is carried out according to the principle: "first is the cheapest transportation, after more expensive", through warehouses $x_8$ and $x_{11}$:

| | **Required amount $\beta_1$** | **From the factory $X_7$** | **From the factory $X_{18}$** | **Required amount $\beta_2$** | **From the factory $X_7$** | **From the factory $X_{18}$** | **Required amount $\beta_3$** | **From the factory $X_7$** | **From the factory $X_{18}$** |
|---|---|---|---|---|---|---|---|---|---|
| $X_{14}$ | 5 | 5 | 0 | 4 | 4 | 0 | 3 | 3 | 0 |
| $X_{15}$ | 5 | 2 | 3 | 3 | 3 | 0 | 4 | 4 | 0 |
| $X_{16}$ | 4 | 0 | 4 | 5 | 3 | 2 | 4 | 0 | 4 |
| $X_{17}$ | 3 | 0 | 3 | 5 | 4 | 1 | 5 | 0 | 5 |
| **Total shipping costs** | | | | | | | | | |
| $X_{14}$ | | 10 | 0 | | 16 | 0 | | 9 | 0 |
| $X_{15}$ | | 6 | 15 | | 18 | 0 | | 20 | 0 |
| $X_{16}$ | | 0 | 16 | | 24 | 16 | | 0 | 24 |
| $X_{17}$ | | 0 | 12 | | 40 | 10 | | 0 | 40 |
| **Total** | | | | | | | | | 276 |

Income from storage in warehouses:

| | **Cost of purchase at the factory $X_7$** | **Amount** | **Total** | **Cost of purchase at the factory $X_{12}$** | **Amount** | **Total** |
|---|---|---|---|---|---|---|
| $\beta_1$ | 38,15 | 7,00 | 267,04 | 45,41 | 10,00 | 454,15 |
| $\beta_2$ | 67,87 | 10,00 | 678,65 | 64,78 | 7,00 | 453,46 |
| $\beta_3$ | 62,29 | 10,00 | 622,87 | 58,22 | 6,00 | 349,35 |
| Total | | | | | | 2825,52 |
| Netprofit | | | | | | 670,90 |

The net profit of the first agent will be **1838.70**

4) ($X_{12}$,$X_{13}$)

Production will be distributed as follows:

| **Production revenue** | **At the $X_{12}$factory** | **Net profit at the factory $X_{12}$** | **At the $X_{13}$factory** | **Net profit at the factory $X_{13}$** | **Total:** |
|---|---|---|---|---|---|
| $\beta_1$ | 10 | 47,82 | 7 | 20,75 | 68,57 |
| $\beta_2$ | 10 | 88,65 | 7 | 18,87 | 107,52 |
| $\beta_3$ | 6 | 29,23 | 10 | 102,87 | 132,10 |
| Total: | 26 | 165.7 | 24 | 142.49 | 308,19 |

Net profit of stores will be:

| | The cost of the factory $X_{12}$ | Quantity of goods | Storage costs | Total Cost Per Unit | Total | The cost of the factory $X_{13}$ | Quantity ofgoods | Storage costs | Total Cost Per Unit | Total | **Net profit** |
|---|---|---|---|---|---|---|---|---|---|---|---|
| $\beta_1$ | 41,78 | 10,00 | 19,00 | 60,78 | 607,82 | 39,96 | 7,00 | 19,00 | 58,96 | 412,75 | 424,43 |
| $\beta_2$ | 67,87 | 10,00 | 33,00 | 100,87 | 1008,65 | 61,70 | 7,00 | 33,00 | 94,70 | 662,87 | 453,48 |
| $\beta_3$ | 56,87 | 6,00 | 22,00 | 78,87 | 473,23 | 62,29 | 10,00 | 22,00 | 84,29 | 842,87 | 523,90 |
| Total | | | | | | | | | | 4008,19 | 1401,81 |

Consider the required amount of resources for each production:

| Factory $x_{12}$ | | | | | | Factory $x_{13}$ | | | | | |
|---|---|---|---|---|---|---|---|---|---|---|---|
| $\beta_1$ | | $\beta_2$ | | $\beta_3$ | | $\beta_1$ | | $\beta_2$ | | $\beta_3$ | |
| 10 | | 10 | | 6 | | 7 | | 7 | | 10 | |
| The required amount of resource per unit of production | | | | | | | | | | | |
| $\alpha_1$ | $\alpha_2$ | $\alpha_1$ | $\alpha_2$ | $\alpha_1$ | $\alpha_2$ | $\alpha_1$ | $\alpha_2$ | $\alpha_1$ | $\alpha_2$ | $\alpha_1$ | $\alpha_2$ |
| 1 | 1 | 1 | 2 | 2 | 1 | 1 | 1 | 1 | 2 | 2 | 1 |
| The required amount of resource for the entire quantity of products | | | | | | | | | | | |
| 10 | 10 | 10 | 20 | 12 | 6 | 7 | 7 | 7 | 14 | 20 | 10 |
| Total resources $\alpha_1$ | | | Total resources $\alpha_2$ | | | Total resources $\alpha_1$ | | | Total resources $\alpha_2$ | | |
| 32 | | | 36 | | | 34 | | | 31 | | |

The cost of transporting a unit of production from the points of extraction of raw materials to the points of production:

| $x_2$ | | $x_3$ | | $x_4$ | | $x_5$ | |
|---|---|---|---|---|---|---|---|
| $x_{12}$ | $x_{13}$ | $x_{12}$ | $x_{13}$ | $x_{12}$ | $x_{13}$ | $x_{12}$ | $x_{13}$ |
| 7 | 6 | 7 | 8 | 7 | 10 | 7 | 12 |
| 17 | 12 | 13 | 12 | 9 | 12 | 5 | 12 |

It is seen that it is most advantageous to transport through points $x_3$ and $x_5$. Then the profit minus the costs of transportation and storage of raw materials will be:

| | Factory $x_{12}$ | | | | Factory $x_{13}$ | | | |
|---|---|---|---|---|---|---|---|---|
| | **cost per unit of transportation** | **storage costs per unit** | **total units** | **total** | **cost per unit of transportation** | **storage costs per unit** | **total units** | **total** |
| $\alpha_1$ | 7,00 | 0,20 | 32,00 | 236,80 | 6,00 | 0,20 | 34,00 | 217,60 |
| $\alpha_2$ | 5,00 | 0,20 | 36,00 | 194,40 | 12,00 | 0,20 | 31,00 | 384,40 |
| Net profit | | | | | | | | 1430,80 |

The cost of transporting a unit of manufactured products from production points to stores where the sale of final products is carried out, is carried out according to the principle: "first is the cheapest transportation, after more expensive", through warehouses $x_8$ and $x_{11}$:

| | **Required amount $\beta_1$** | **From the factory $X_{12}$** | **From the factory $X_{13}$** | **Required amount $\beta_2$** | **From the factory $X_{12}$** | **From the factory $X_{13}$** | **Required amount $\beta_3$** | **From the factory $X_{12}$** | **From the factory $X_{13}$** |
|---|---|---|---|---|---|---|---|---|---|
| $X_{14}$ | 5 | 0 | 5 | 4 | 0 | 4 | 3 | 0 | 3 |
| $X_{15}$ | 5 | 3 | 2 | 3 | 0 | 3 | 4 | 0 | 4 |
| $X_{16}$ | 4 | 4 | 0 | 5 | 5 | 0 | 4 | 1 | 3 |
| $X_{17}$ | 3 | 3 | 0 | 5 | 5 | 0 | 5 | 5 | 0 |
| **Total shipping costs** | | | | | | | | | |
| $X_{14}$ | | 0 | 15 | | 0 | 24 | | 0 | 12 |
| $X_{15}$ | | 12 | 8 | | 0 | 24 | | 0 | 24 |
| $X_{16}$ | | 12 | 0 | | 30 | 0 | | 5 | 24 |
| $X_{17}$ | | 0 | 0 | | 0 | 0 | | 0 | 0 |
| **Total** | | | | | | | | | 190 |

Income from storage in warehouses:

| | **Cost of purchase at the factory $X_{12}$** | **Amount** | **Total** | **Cost of purchase at the factory $X_{13}$** | **Amount** | **Total** |
|---|---|---|---|---|---|---|
| $\beta_1$ | 41,78 | 10,00 | 417,82 | 39,96 | 7,00 | 279,75 |
| $\beta_2$ | 67,87 | 10,00 | 678,65 | 61,70 | 7,00 | 431,87 |
| $\beta_3$ | 56,87 | 6,00 | 341,23 | 62,29 | 10,00 | 622,87 |
| Total | | | | | | 2772,19 |
| Net profit | | | | | | 681,56 |

The net profit of the first agent will be**1922,3**

5)($X_{12}$,$X_{18}$)

Production will be distributed as follows:

| **Production revenue** | **At the $X_{12}$factory** | **Net profit at the factory $X_{12}$** | **At the $X_{18}$factory** | **Net profit at the factory $X_{18}$** | **Total:** |
|---|---|---|---|---|---|
| $\beta_1$ | 7 | 33,47 | 10 | 84,15 | 117,62 |
| $\beta_2$ | 10 | 88,65 | 7 | 40,46 | 129,11 |
| $\beta_3$ | 6 | 29,23 | 10 | 62,25 | 91,47 |
| Total: | 23 | 151.35 | 27 | 186.86 | 338,21 |

Net profit of stores will be:

| | The cost of the factory $X_{12}$ | Quantity of goods | Storage costs | Total Cost Per Unit | Total | The cost of the factory $X_{18}$ | Quantity of goods | Storage costs | Total Cost PerUnit | Total | **Net profit** |
|---|---|---|---|---|---|---|---|---|---|---|---|
| $\beta_1$ | 41,78 | 7,00 | 19,00 | 60,78 | 425,47 | 45,41 | 10,00 | 19,00 | 64,41 | 644,15 | 375,38 |
| $\beta_2$ | 67,87 | 10,00 | 33,00 | 100,87 | 1008,65 | 64,78 | 7,00 | 33,00 | 97,78 | 684,46 | 431,89 |
| $\beta_3$ | 56,87 | 6,00 | 22,00 | 78,87 | 473,23 | 58,22 | 10,00 | 22,00 | 80,22 | 802,25 | 564,53 |
| Total | | | | | | | | | | 4038,21 | 1371,79 |

Consider the required amount of resources for each production:

| Factory $x_{12}$ | | | | | | Factory $x_{18}$ | | | | | |
|---|---|---|---|---|---|---|---|---|---|---|---|
| $\beta_1$ | | $\beta_2$ | | $\beta_3$ | | $\beta_1$ | | $\beta_2$ | | $\beta_3$ | |
| 7 | | 10 | | 6 | | 10 | | 7 | | 10 | |
| The required amount of resource per unit of production | | | | | | | | | | | |
| $\alpha_1$ | $\alpha_2$ | $\alpha_1$ | $\alpha_2$ | $\alpha_1$ | $\alpha_2$ | $\alpha_1$ | $\alpha_2$ | $\alpha_1$ | $\alpha_2$ | $\alpha_1$ | $\alpha_2$ |
| 1 | 1 | 1 | 2 | 2 | 1 | 1 | 1 | 1 | 2 | 2 | 1 |
| The required amount of resource for the entire quantity of products | | | | | | | | | | | |
| 7 | 7 | 10 | 20 | 12 | 6 | 10 | 10 | 7 | 14 | 20 | 10 |
| Total resources $\alpha_1$ | | | Total resources $\alpha_2$ | | | Total resources $\alpha_1$ | | | Total resources $\alpha_2$ | | |
| 29 | | | 33 | | | 37 | | | 34 | | |

The cost of transporting a unit of production from the points of extraction of raw materials to the points of production:

| | $x_2$ | | $x_3$ | | $x_4$ | | $x_5$ | |
|---|---|---|---|---|---|---|---|---|
| | $x_{12}$ | $x_{18}$ | $x_{12}$ | $x_{18}$ | $x_{12}$ | $x_{18}$ | $x_{12}$ | $x_{18}$ |
| $\alpha_1$ | 7 | 8 | 7 | 8 | 7 | 8 | 7 | 8 |
| $\alpha_2$ | 17 | 16 | 13 | 14 | 9 | 10 | 5 | 6 |

It is seen that it is most advantageous to transport through points $x_4$ and $x_5$. Then the profit minus the costs of transportation and storage of raw materials will be:

| | Factory $x_{12}$ | | | | Factory $x_{18}$ | | | |
|---|---|---|---|---|---|---|---|---|
| | **cost per unit of transportation** | **storage costs per unit** | **total units** | **total** | **cost per unit of transportation** | **storage costs per unit** | **total units** | **total** |
| $\alpha_1$ | 7,00 | 0,20 | 29,00 | 214,60 | 8,00 | 0,20 | 37,00 | 310,80 |
| $\alpha_2$ | 9,00 | 0,20 | 33,00 | 310,20 | 6,00 | 0,20 | 34,00 | 217,60 |
| Net profit | | | | | | | | 1410,80 |

The cost of transporting a unit of manufactured products from production points to stores where the sale of manufactured productsis carried out, is carried out according to the principle: “first, the cheapest transportation, after more expensive”, through warehouses $x_{10}$ and $x_{11}$:

| | **Required amount $\beta_1$** | **From the factory $X_{12}$** | **From the factory $X_{18}$** | **Required amount $\beta_2$** | **From the factory $X_{12}$** | **From the factory $X_{18}$** | **Required amount $\beta_3$** | **From the factory $X_{12}$** | **From the factory $X_{18}$** |
|---|---|---|---|---|---|---|---|---|---|
| $X_{14}$ | 5 | 5 | 0 | 4 | 4 | 0 | 3 | 3 | 0 |
| $X_{15}$ | 5 | 2 | 3 | 3 | 3 | 0 | 4 | 3 | 1 |
| $X_{16}$ | 4 | 0 | 4 | 5 | 3 | 2 | 4 | 0 | 4 |
| $X_{17}$ | 3 | 0 | 3 | 5 | 0 | 5 | 5 | 0 | 5 |
| **Total shipping costs** | | | | | | | | | |
| $X_{14}$ | | 25 | 0 | | 40 | 0 | | 27 | 0 |
| $X_{15}$ | | 8 | 15 | | 24 | 0 | | 21 | 8 |
| $X_{16}$ | | 0 | 12 | | 18 | 16 | | 0 | 24 |
| $X_{17}$ | | 0 | 12 | | 0 | 50 | | 0 | 40 |
| **Total** | | | | | | | | | 340 |

Income from storage in warehouses:

| | **Cost of purchase at the factory $X_{12}$** | **Amount** | **Total** | **Cost of purchase at the factory $X_{18}$** | **Amount** | **Total** |
|---|---|---|---|---|---|---|
| $\beta_1$ | 41,78 | 7,00 | 292,47 | 45,41 | 10,00 | 454,15 |
| $\beta_2$ | 67,87 | 10,00 | 678,65 | 64,78 | 7,00 | 453,46 |
| $\beta_3$ | 56,87 | 6,00 | 341,23 | 58,22 | 10,00 | 582,25 |
| Total | | | | | | 2802,21 |
| Net profit | | | | | | 675,56 |

The net profit of the first agent will be **1746.36**

6)($X_{13}$,$X_{18}$)

Production will be distributed as follows:

| **Production revenue** | **At the $X_{13}$factory** | **Net profit at the factory $X_{13}$** | **At the $X_{18}$factory** | **Net profit at the factory $X_{18}$** | **Total:** |
|---|---|---|---|---|---|
| $\beta_1$ | 7 | 20,75 | 10 | 84,15 | 104,90 |
| $\beta_2$ | 7 | 18,87 | 10 | 57,80 | 76,67 |
| $\beta_3$ | 10 | 102,87 | 6 | 37,35 | 140,22 |
| Total: | 24 | 142.49 | 26 | 179.3 | 321,79 |

Net profit of stores will be:

| | The cost of the factory $X_{13}$ | Quantity of goods | Storage costs | Total Cost Per Unit | Total | The cost of the factory $X_{18}$ | Quantity of goods | Storage costs | Total Cost Per Unit | Total | **Net profit** |
|---|---|---|---|---|---|---|---|---|---|---|---|
| $\beta_1$ | 39,96 | 7,00 | 19,00 | 58,96 | 412,75 | 45,41 | 10,00 | 19,00 | 64,41 | 644,15 | 388,10 |
| $\beta_2$ | 61,70 | 7,00 | 33,00 | 94,70 | 662,87 | 64,78 | 10,00 | 33,00 | 97,78 | 977,80 | 484,33 |
| $\beta_3$ | 62,29 | 10,00 | 22,00 | 84,29 | 842,87 | 58,22 | 6,00 | 22,00 | 80,22 | 481,35 | 515,78 |
| Total | | | | | | | | | | 4021,79 | 1388,21 |

Consider the required amount of resources for each production:

| Factory $x_{13}$ | | | | | | Factory $x_{18}$ | | | | | |
|---|---|---|---|---|---|---|---|---|---|---|---|
| $\beta_1$ | | $\beta_2$ | | $\beta_3$ | | $\beta_1$ | | $\beta_2$ | | $\beta_3$ | |
| 7 | | 7 | | 10 | | 10 | | 10 | | 6 | |
| The required amount of resource per unit of production | | | | | | | | | | | |
| $\alpha_1$ | $\alpha_2$ | $\alpha_1$ | $\alpha_2$ | $\alpha_1$ | $\alpha_2$ | $\alpha_1$ | $\alpha_2$ | $\alpha_1$ | $\alpha_2$ | $\alpha_1$ | $\alpha_2$ |
| 1 | 1 | 1 | 2 | 2 | 1 | 1 | 1 | 1 | 2 | 2 | 1 |

| The required amount of resource for the entire quantity of products | | | | | | | | | | | |
|---|---|---|---|---|---|---|---|---|---|---|---|
| 7 | 7 | 7 | 14 | 20 | 10 | 10 | 10 | 10 | 20 | 12 | 6 |
| Total resources $\alpha_1$ | | | Total resources $\alpha_2$ | | | Total resources $\alpha_1$ | | | Total resources $\alpha_2$ | | |
| 34 | | | 31 | | | 32 | | | 36 | | |

The cost of transporting a unit of production from the points of extraction of raw materials to the points of production:

| | $x_2$ | | $x_3$ | | $x_4$ | | $x_5$ | |
|---|---|---|---|---|---|---|---|---|
| | $x_{12}$ | $x_{18}$ | $x_{12}$ | $x_{18}$ | $x_{12}$ | $x_{18}$ | $x_{13}$ | $x_{18}$ |
| $\alpha_1$ | 6 | 8 | 8 | 8 | 10 | 8 | 12 | 8 |
| $\alpha_2$ | 12 | 16 | 12 | 14 | 12 | 10 | 12 | 6 |

It is seen that it is most advantageous to transport through points $x_2$ and $x_5$. Then the profit minus the costs of transportation and storage of raw materials will be:

| | Factory $x_{13}$ | | | | Factory $x_{18}$ | | | |
|---|---|---|---|---|---|---|---|---|
| | **cost per unit of transportation** | **storage costs per unit** | **total units** | **total** | **cost per unit of transportation** | **storage costs per unit** | **total units** | **total** |
| $\alpha_1$ | 6,00 | 0,20 | 34,00 | 217,60 | 8,00 | 0,20 | 32,00 | 268,80 |
| $\alpha_2$ | 12,00 | 0,20 | 31,00 | 384,40 | 6,00 | 0,20 | 36,00 | 230,40 |
| Net profit | | | | | | | | 1362,80 |

The cost of transporting a unit of manufactured products from production points to stores where the sale of final products is carried out, is carried out according to the principle: “first is the cheapest transportation, after more expensive”, through warehouses $x_8$ and $x_{11}$:

| | **Required amount $\beta_1$** | **From the factory $X_{13}$** | **From the factory $X_{18}$** | **Required amount $\beta_2$** | **From the factory $X_{13}$** | **From the factory $X_{18}$** | **Required amount $\beta_3$** | **From the factory $X_{13}$** | **From the factory $X_{18}$** |
|---|---|---|---|---|---|---|---|---|---|
| $X_{14}$ | 5 | 5 | 0 | 4 | 4 | 0 | 3 | 3 | 0 |
| $X_{15}$ | 5 | 2 | 3 | 3 | 3 | 0 | 4 | 4 | 0 |
| $X_{16}$ | 4 | 0 | 4 | 5 | 0 | 5 | 4 | 3 | 1 |
| $X_{17}$ | 3 | 0 | 3 | 5 | 0 | 5 | 5 | 0 | 5 |
| **Total shipping costs** | | | | | | | | | |
| $X_{14}$ | | 45 | 0 | | 72 | 0 | | 48 | 0 |
| $X_{15}$ | | 16 | 15 | | 48 | 0 | | 56 | 0 |
| $X_{16}$ | | 0 | 20 | | 0 | 40 | | 36 | 6 |
| $X_{17}$ | | 0 | 12 | | 0 | 50 | | 0 | 40 |
| **Total** | | | | | | | | | 504 |

Income from storage in warehouses:

| | **Cost of purchase at the factory $X_{12}$** | **Amount** | **Total** | **Cost of purchase at the factory $X_{18}$** | **Amount** | **Total** |
|---|---|---|---|---|---|---|
| $\beta_1$ | 39,96 | 7,00 | 279,75 | 45,41 | 10,00 | 454,15 |
| $\beta_2$ | 61,70 | 7,00 | 431,87 | 64,78 | 10,00 | 647,80 |
| $\beta_3$ | 62,29 | 10,00 | 622,87 | 58,22 | 6,00 | 349,35 |
| Total | | | | | | 2785,79 |
| Net profit | | | | | | 678,84 |

The net profit of the first agent will be 1537.64

Knowing the income of each agent, we can calculate a compromise solution for their distribution.

**Step 1**. Build a matrix of agent revenues depending on their choice of location:

| | (7,12) | (7,13) | (7,18) | (12,13) | (12,18) | (13,18) |
|---|---|---|---|---|---|---|
| First agent | 1963,47 | 1654,04 | 1838,70 | 1922,36 | 1746,36 | 1537,64 |
| Second agent | 338,66 | 309,80 | 361,52 | 308,19 | 338,21 | 321,79 |
| Third agent | 1371,34 | 1400,20 | 1348,81 | 1401,81 | 1371,79 | 1388,21 |

**Step 2.** Create an "ideal vector", consisting of the maximum income received by sellers.

$$M = \begin{pmatrix} M_1 \\ \vdots \\ M_l \end{pmatrix}, M_l = \max_m \alpha_{l,m}$$

$$M = \begin{vmatrix} M_1 \\ M_2 \\ M_3 \end{vmatrix} = \begin{vmatrix} 1963{,}47 \\ 361{,}52 \\ 1401{,}81 \end{vmatrix}$$

**Step 3.** Calculate the "discrepancy" - the magnitude of the deviations of income from the maximum income for each seller

$$\Gamma_M = (M - \alpha_{l,m}) = (\beta_{l,m})$$

| | (7,12) | (7,13) | (7,18) | (12,13) | (12,18) | (13,18) |
|---|---|---|---|---|---|---|
| First agent | 0 | 309,43 | 124,77 | 41,11 | 217,11 | 425,83 |
| Second agent | 22,86 | 51,72 | 0 | 53,33 | 23,31 | 39,73 |
| Third agent | 30,47 | 1,61 | 53 | 0 | 30,02 | 13,60 |

**Step 4.** In each situation, we will sort the incomes in ascending order in such a way that the first line will contain the smallest amounts of income, and the bottom line will contain the largest ones. Thus, the last line will be:

$$\max_l\left(\beta_{m,l}\right)=\max_l\left(M-\alpha_{l,m}\right)$$

| | (7,12) | (7,13) | (7,18) | (12,13) | (12,18) | (13,18) |
|---|---|---|---|---|---|---|
| | 0 | 1,61 | 0 | 0 | 23,31 | 13,6 |
| | 22,86 | 51,72 | 53 | 53,33 | 30,02 | 39,73 |
| | 30,47 | 309,43 | 124,77 | 41,11 | 217,11 | 425,83 |

**Step 5.** Among the maximum residuals found, select the minimum value.

If there are several situations in the last line with the minimum found, then go up one line and there we are already looking for a minimum, and so on. The situations thus obtained (the situation) will be a completely compromise set or a compromise solution:

| (13,18) |
|---|
| 13,6 |
| 39,73 |
| 425,83 |

The winnings of the agents are as follows: (1537.64; 321.79; 1388.21).The location of production points at points $x_{13}$, $x_{18}$, points of temporary storage of raw materials at points $x_2$, $x_5$, temporary storage of final products at points $x_8$, $x_{11}$.

## 9. Conclusion

Thus, the proposed method for solving optimization problems related to the territorial location of warehouses, raw materials processing plants and shops, agent costs, production planning can be used in real life and bring good profits to certain agents / companies.

## Acknowledgements

The work is partly supported by work RFBR No. 18-01-00796.